\documentstyle[11pt]{article}

\input{psfig}

\begin{document}
\newcommand{\bh}{the boundary of hyperbolicity}
\newcommand{\hs}{hyperbolic space}
\newcommand{\ds}{dynamical system}
\newcommand{\dss}{dynamical systems}
\newcommand{\cp}{critical point}
\newcommand{\nl}{nonlinear}
\newcommand{\tr}{transformation}
\newcommand{\trs}{transformations}
\newcommand{\fn}{function}
\newcommand{\fns}{functions}
\newcommand{\acim}{absolutely continuous invariant measure}
\newcommand{\intl}{the interval $[-1, 1]$}
\newcommand{\pl}{power law}
\newcommand{\sm}{\mbox{$dy
=dx/((1+\varepsilon)^{2}-x^{2})^{\frac{\gamma-1}{\gamma}}$}}
\newcommand{\sms}{\mbox{$dx/((1+\varepsilon)^{2}-x^{2})^{\frac{\gamma-1}{\gamma}}$}}
\newcommand{\st}{\mbox{$w=i_{0}i_{1}\cdots i_{n}$}}
\newcommand{\npt}{\mbox{$n^{th}$}-partition}
\newcommand{\pts}{partitions}
\newcommand{\ep}{\varepsilon}
\newcommand{\al}{\mbox{$C^{1+\alpha }$}}
\newcommand{\e}{embedding}
\newcommand{\ef}{exponentially fast}
\newcommand{\mliF}{
the maximal length of the intervals in
the \mbox{$n^{th}$}-partition \mbox{$\eta_{n, F}$} goes to zero}
\newcommand{\mli}{
the maximal length of intervals in
\mbox{$\eta_{n, f_{\varepsilon}}$} goes to zero}
\newcommand{\ngi}{as n goes to infinite}
\newcommand{\sd}{Schwarzian derivative}
\newcommand{\nsd}{negative Schwarzian derivative}
\newcommand{\psd}{positive Schwarzian derivative}
\newcommand{\stp}{sequence of partitions \mbox{$\{
\eta_{n,f}\}_{n=0}^{\infty}$}}
\newcommand{\bg}{bounded geometry}
\newcommand{\ode}{ordinary differential equation}
\newcommand{\odi}{ordinary differential inequality}
\newcommand{\nd}{\noindent}
\newcommand{\T}{Theorem}
\newcommand{\Le}{Lemma}
\newcommand{\C}{Corollary}
\newcommand{\ex}{expanding}
\newcommand{\poc}{positive constant}
\newcommand{\pocs}{positive constants}
\newcommand{\poi}{positive integer}
\newcommand{\pis}{positive integers}
\newcommand{\mi}{the minimal value}
\newcommand{\tc}{condition}
\newcommand{\tcs}{conditions}
\newcommand{\pn}{positive number}
\newcommand{\pns}{positive numbers}
\newcommand{\Q}{{\sl QED.}}

\vspace*{1in}
\centerline {\Large \bf Dynamics of Certain Smooth One-dimensional 
Mappings}

\vskip10pt
\centerline {\Large \bf IV. Asymptotic geometry of Cantor sets}

\vskip20pt
\large
\centerline{Yunping Jiang }
\centerline{Institute for Mathematical Sciences, SUNY at Stony Brook}
\centerline{Stony Brook, L.I., NY 11794}
\vskip20pt
\centerline{May 4, 1991}

\vskip20pt
\centerline{\em Dedicated to Professor John Milnor on the occasion of his sixtieth birthday}

\vskip50pt
\centerline{\Large \bf Abstract}

\vskip20pt
We study hyperbolic mappings
depending on a parameter $\varepsilon $. Each of them has
an invariant Cantor set.
As $\varepsilon $
tends to
zero, the mapping approaches the boundary of
hyperbolicity.
We
analyze the
asymptotics of the gap geometry and the scaling function
geometry of the invariant Cantor
set as
$\varepsilon $ goes to zero. For example, in the quadratic case, we
show that all the gaps close uniformly with speed $\sqrt {\varepsilon
}$. There is a limiting scaling function of the limiting mapping
and this scaling function has dense jump discontinuities because
the limiting
mapping is not expanding.  Removing these discontinuities by
continuous extension, we show that we obtain the
scaling function of the limiting mapping with respect to the
Ulam-von Neumann type metric.

\pagebreak

\vspace*{1in}

\centerline{\Large \bf Contents}

\vskip30pt
\noindent{\S 1 Introduction}

\vskip10pt
\noindent{\S 2 The Boundary of Hyperbolicity, ${\cal BH}$}
 
\vskip5pt
\S 2.1 The singular change of metric on the interval

\vskip5pt
\S 2.2 The definition of the boundary of hyperbolicity

\vskip10pt
\noindent{\S 3 The Space of Hyperbolic Mappings, ${\cal H}$}

\vskip5pt
\S 3.1 The definition of the scaling function 

\vskip5pt
\S 2.4 The definition of the space of hyperbolic mappings

\vskip10pt
\noindent{\S 3 Asymptotic Geometry of Cantor Sets}
 
\vskip5pt
\S 3.1 Good families of mappings in ${\cal BH}\cup {\cal H}$

\vskip5pt
\S 3.2 Asymptotic scaling function geometry of Cantor sets
 
\vskip5pt
\S 3.3 Determination of the geometry of Cantor set by the leading \nolinebreak gap

\vskip5pt
\S 3.4 Proof of Theorem A

\vskip5pt
\S 3.5 Scaling functions of mappings on ${\cal BH}$

\vskip5pt
\S 3.6 The Hausdorff dimension of the maximal invariant set of $q_{\varepsilon}$

\pagebreak

\noindent{\Large \bf \S 1 Introduction }

\vskip10pt
Ulam and von Neumann studied the nonlinear self mapping $q(x)=1-2x^{2}$ of the interval $[-1, 1]$. They observed that
$\rho_{q}=1/(\pi \sqrt {1-x^{2}})$ is
the density
function of a unique absolutely continuous $q$-invariant measure (we
only consider probability measures). In modern language,
this observation follows from making
the singular change of metric $|dy| = (2|dx|)/(\pi \sqrt {1-x^{2}})$.
If we let $y=h(x)$ be the corresponding change of coordinate and
$\tilde{q}=h\circ
q\circ h^{-1}$, then $q$ becomes $\tilde{q}(y)=1-2|y|$,
a piecewise linear mapping with expansion rate $2$ on $[-1,
1]$.  The dynamics of $\tilde{q}$ is more easily understood.

\vskip5pt
\centerline{\psfig{figure=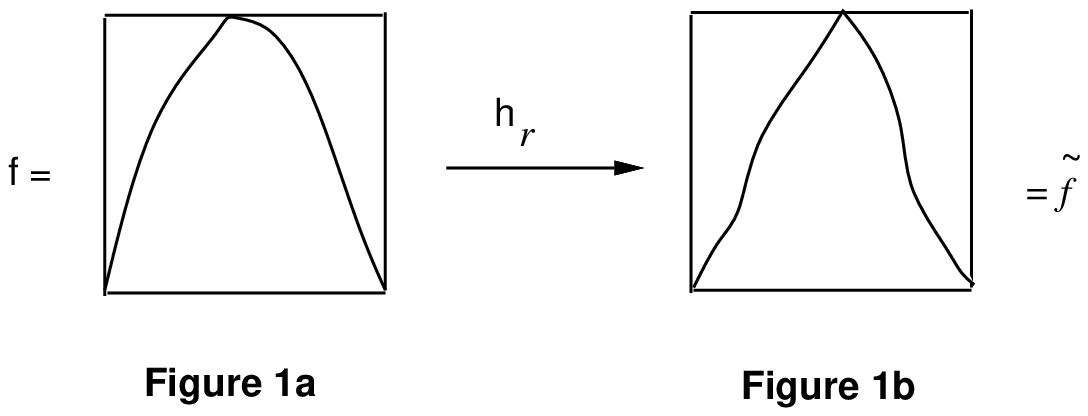}}

\vskip5pt
Under quite general conditions, a mapping whose graph
looks like the graph shown in Figure 2 has hyperbolic properties.
We may say that a mapping whose graph looks like the graph shown in
Figure 1a is on the \underline {boundary of hyperbolicity} (a more
precise definition of what we call the boundary of hyperbolicity is
given below).

\vskip5pt
\centerline{\psfig{figure=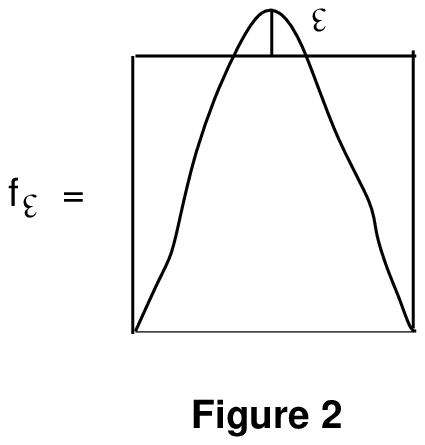}}

\vskip5pt
In order to study
more general smooth self mappings
of the interval with a unique power law critical point,
we employ a
change of metric similar to the one used by Ulam and von Neumann.
The \underline {change of}
\underline {metric has singularities
 of \hspace{.05pt} same type \hspace{.05pt} at the two boundary
\hspace{.05pt} points of}

\noindent \underline {the interval}.
It is universal in the sense that it does not depend on
particular mapping
$f$, but only on the power law $|x|^{\gamma}$ at the critical point.
Suppose $y=h_{\gamma}(x)$ is the corresponding change of coordinate on
the interval.  After this change of coordinate,
$f$ becomes $\tilde {f}=h_{\gamma}\circ f\circ h_{\gamma }^{-1}$ (see
Figure $1b$), which is smooth
except at the critical point.  The mapping $\tilde{f}$ has nonzero
derivative at every point except the critical point.  At
the critical point,  the left and the right derivatives
of $\tilde{f}$ exist and are positive and negative, respectively.

\vskip5pt
A nice feature of the mapping $q(x)=1-2x^{2}$
is that $\tilde{q}$ is expanding
with H\"older continuous derivative,
which implies that a
certain
binary tree of intervals associated with the dynamics of $\tilde{q}$
has bounded geometry (see [J2]).  The expanding property does not
carry
over to our more general setting but the bounded geometry does (see [J2]).

\vskip5pt
Suppose $E_{0}$
is the set consisting of the critical point of $f$ and of the
two boundary points of the
interval.  For every positive integer $n$,  let $E_{n}$ be the preimage
of $E_{n-1}$ under $f$.
The $n^{th}$-partition $\eta_{n}$ of the interval determined by
$f$ is the collection of all the subintervals
bounded by consecutive points
of $E_{n}$.  Let $\lambda_{n}$ be the maximum length of the intervals in
the $n^{th}$-partition.

\vskip5pt
We say the sequence of nested partitions $\{ \eta_{n}\}_{n=0}^{\infty
}$ determined by $f$
decreases exponentially if $\lambda_{n}$
decreases exponentially.

\vskip5pt
We need the following conditions in this paper.

\vskip5pt
\noindent (*) $f$ to be a $C^{1}$ self mapping
of an interval with a unique power law, $|x|^{\gamma}$ for some
$\gamma >1$, critical point and to map
the critical point to the right endpoint of the interval and both
endpoints of the interval to the left endpoint (see Figure 1a),

\vskip5pt
\noindent (**) the derivative of $\tilde{f}$ to be
piecewise $\alpha $-H\"older continuous for some
$0< \alpha \leq 1$
and

\vskip5pt
\noindent (***) the sequence of nested partitions $\{
\eta_{n}\}_{n=0}^{\infty}$ determined by $f$ to decrease exponentially.

\vskip5pt
The exponential decay of (***) were proved in [J2] under
either of the following two hypotheses.

\vskip5pt
{\sl (1)} The mapping $f$ is $C^{3}$ with nonpositive
Schwarzian derivative and expanding at both
boundary points of the interval, that is, the absolute values of
the derivatives of
$f$ at both boundary points are greater than one (see also
[Mi]).
The Schwarzian derivative of $f$ is $S(f)=f'''/f'-(3/2)(f''/f')^{2}$.

\vskip5pt
{\sl (2)} The mapping $\tilde {f}^{\prime }$ is piecewise
Lipschitz
and all the periodic points of $f$ are expanding, that is, the
absolute values of the eigenvalues of $f$ at all periodic points are
greater than
one (see also [Ma]). The eigenvalue of $f$ at a periodic point $p$ of
period $n$ of $f$ is $e_{p} =(f^{\circ n})'(p)$.

\vskip5pt
We call the set of mappings
satisfying (*), (**) and (***) the
\underline {boundary of hyperbolicity, ${\cal BH}$}.	

\vskip20pt
The
mappings $f$ on ${\cal BH}$ are limits of
mappings $f_{\varepsilon }$,  which do not keep the
interval invariant,  but keep invariant a Cantor
set $\Lambda_{\varepsilon }$ having bounded geometry (see
Figure 2). We say the mappings like $f_{\varepsilon }$ are
hyperbolic
(a more precise definition
of what we call a hyperbolic mapping is given
below).
The \underline {space of hyperbolic mappings}, as well as
the asymptotic behavior of these hyperbolic mappings as they approach
the boundary of hyperbolicity, is the topic of this paper.

\vskip5pt
We use $f_{\varepsilon , 0}$ and
$f_{\varepsilon , 1}$ to denote the left and right branches of
$f_{\varepsilon }$, respectively.
Let $g_{\varepsilon, 0}$
and $g_{\varepsilon, 1}$ be the inverses of
$f_{\varepsilon , 0}$ and $f_{\varepsilon ,1}$.
For a finite string $w=i_{0}\cdots i_{n}$ of zeroes and ones,
we use $g_{\varepsilon , w}$ to denote the
composition
$g_{\varepsilon , w} =
g_{\varepsilon , i_{0}}\circ
\cdots \circ g_{\varepsilon , i_{n}}$.
Define $I_{\varepsilon , w}$ to be the image under $g_{\varepsilon , w}$ of
the interval where $f_{\ep}$ is defined
(see Figure 3).

\vskip5pt
\begin{picture}(300,50)(0,0)
\put(20,40){\line(1,0){100}}
\put(65,42){$I_{\varepsilon , w}$}
\put(20,25){\line(1,0){40}}
\put(35,27){$I_{\ep ,w0}$}
\put(80,25){\line(1,0){40}}
\put(95,27){$I_{\ep , w1}$}
\put(200,40){\line(1,0){100}}
\put(245,42){$I_{\ep , w}$}
\put(200,25){\line(1,0){40}}
\put(215,27){$I_{\ep , w1}$}
\put(260,25){\line(1,0){40}}
\put(275,27){$I_{\ep , w0}$}
\put(150,35){or}
\put(130,0){Figure 3}
\end{picture}

\vskip5pt
Suppose $\eta_{n,
\varepsilon }$ is the
collection of $I_{\varepsilon , w}$ for all finite strings $w$ of
zeroes and ones of length $n+1$.  We use $\lambda_{n,\varepsilon }$
to denote the maximum length of the intervals in $\eta_{n,\varepsilon
}$. Notice that the union of all the intervals in
$\eta_{n, \varepsilon }$ covers the maximal invariant set
of $f_{\varepsilon }$.

\vskip5pt
\begin{picture}(300,70)(0,0)
\put(20,60){\line(1,0){100}}
\put(65,62){$.w$}
\put(20,45){\line(1,0){40}}
\put(35,47){$.w0$}
\put(80,45){\line(1,0){40}}
\put(95,47){$.w1$}
\put(70, 30){$\vdots$}
\put(200,60){\line(1,0){100}}
\put(245,62){$w.$}
\put(200,45){\line(1,0){40}}
\put(220,47){$w1.$}
\put(260,45){\line(1,0){40}}
\put(275,47){$w0.$}
\put(250,30){$\vdots$}
\put(65, 15){${\cal C}$}
\put(245,15){${\cal C}^{*}$}
\put(50,0){Figure 4a}
\put(230,0){Figure 4b}
\end{picture}

\vskip5pt
There are two
topologies on the set of all the labellings $w$.
One topology is induced by
reading the labellings $w$ from left to right;
the other by reading
the labellings $w$ from right to
left.  The limit set of the set of these labellings $w$ in the topology
induced by reading the labellings $w$ from right to left is
the phase space of the dynamical system
$f_{\varepsilon }$.  We call it the \underline {topological Cantor set}
${\cal C}$ (see Figure 4a).  Points in ${\cal C}$ are
one-sided infinite strings of zeroes and ones extending infinitely
to the right.  If we take the limit set of the
labellings $w$ in the topology induced by reading from right to
left, we obtain the \underline {dual Cantor set} ${\cal C}^{*}$.
A point in ${\cal C}^{*}$ is called a ``dual point'' which is
one-sided infinite string
of zeroes and ones extending infinitely to the left.

\vskip5pt
The \underline {scaling function} of $f_{\ep}$, when it is defined, is a
function defined on
${\cal C}^{*}$.  Assume $a^{*} \in {\cal C}^{*}$, so that $a^{*}$ is a
one-sided infinite string of zeroes and ones extending infinitely to
the left. Suppose
$a^{*}=(\cdots wi.)$, where $w$ is a finite string of zeroes and
ones
and $i$ is either zero or one. Note that $I_{wi}$ is a subinterval of
$I_{w}$.  Let $s(wi)$ equal the ratio of the lengths,
$|I_{wi}|/|I_{w}|$. We let $s(a^{*})$ be the limit set of $s(wi)$ as the
length of $w$ tends to infinity.  If this limit set consists of just one
number for every $a^{*}\in {\cal C}$, then we say that $s(a^{*})$ is the
scale of $f_{\ep}$ at $a^{*}$ and that $s$ is the scaling function of
$f_{\ep}$ defined on ${\cal C}^{*}$.
Note that the scaling function $s(a^{*})$ of $f_{\ep}$ depends on $\ep$.
Sometimes we denote it by $s_{\ep}(a^{*})$.  The same definition gives
the scaling function $s_{0}(a^{*})$ of a mapping $f_{0}$ on ${\cal BH}$.
Since for the mappings $f_{\ep}$, $\ep
\geq 0$, the length of the interval $I_{w}$ converges to zero uniformly
as the length of $w$ approaches infinity, it is obvious that the scaling
function is a $C^{1}$-invariant. Recall that a smooth invariant is an
object associated to $f_{\ep}$ which is the same for $f_{\ep}$ as for
$h\circ f_{\ep}\circ h^{-1}$ whenever $h$ is an orientation preserving
$C^{1}$-diffeomorphism.

\vskip5pt
To be sure that the limits defining the scaling function of
$f_{\ep}$ actually exist, we need

\vskip5pt
\noindent $(i)$ $f_{\ep}$ to be a mapping from an interval to the real
line which maps its unique
critical point
out of this interval and both endpoints of the interval
to the left endpoint (see Figure 2),

\vskip5pt
\noindent $(ii)$ $f_{\varepsilon }$ to be $C^{1+\alpha }$ for
some $0< \alpha \leq 1$ and

\vskip5pt
\noindent $(iii)$
the sequence of maximum lengths
$\lambda_{n, \varepsilon }$ determined by
$f_{\varepsilon }$ to decrease exponentially.

\vskip5pt
Parallel to what we do for the boundary of hyperbolicity, we prove the
exponential decay of $(iii)$ under either of the following two
hypotheses.

\vskip5pt
{\sl (1)} The mapping $f_{\ep}$ is $C^{3}$ with nonpositive
Schwarzian derivative and expanding at both
boundary points of the interval.

\vskip5pt
{\sl (2)} The mapping $f_{\ep}$ is $C^{1,1}$ and all the
periodic points of $f_{\ep}$ are expanding.

\vskip5pt
We call the set of mappings $f_{\ep}$ satisfying $(i)$, $(ii)$ and
$(iii)$ the \underline {space of hyperbolic
mappings, ${\cal H}$}. Without
loss of generality, we may assume that the interval where $f_{\ep}$ is
defined is the interval $[-1,1]$ and that the critical point of
$f_{\ep}$ is zero.

\vskip5pt
Sullivan [S2] showed that
the scaling function is a complete invariant for $C^{1}$-conjugacy of
mappings in ${\cal H}$.
It plays the same role that eigenvalues play in the $C^{2}$-case
(recall that Sullivan [S1] showed that the
eigenvalues are complete invariants for $C^{1}$-conjugacy
of $C^{2}$ mappings of ${\cal H}$).
We examine
the asymptotics of the scaling function of
$f_{\varepsilon}$ as $\varepsilon $ decreases to zero.

\vskip5pt
\begin{picture}(300,125)(0,0)
\put(150,105){\line(0,1){8}}
\put(152,107){$\varepsilon$}
\put(100,15){\framebox(90,90)}
\put(30, 55){$\{ f_{\varepsilon} \}_{0\leq \varepsilon \leq \ep_{0}} =$}
\put(125,0){Figure 5}
\end{picture}

\vskip5pt
Hereafter, we say that a family
$\{ f_{\varepsilon } \}_{0\leq \varepsilon \leq \ep_{0}}$ (see Figure 5)
of mappings in ${\cal BH}\cup {\cal H}$ is a \underline {good family} if
it satisfies the following conditions:

\vskip5pt
\begin{enumerate}
\item the family $f_{\ep }(x)$ is $C^{1}$ in both variables $\ep$ and
$x$,

\item each $ f_{\ep}$ has the same power law
$|x|^{\gamma }$ with $\gamma >1$ at
the critical points and
$R^{-}(x, \ep)
=f_{\ep}'(x)/|x|^{\gamma -1}$ defined on $[-1,0]\times [0,\ep_{0}]$
and $R^{+}(x, \ep)=f_{\ep}'(x)/|x|^{\gamma -1}$ defined on $[0,1]\times
[0,\ep_{0}]$ are continuous,

\item there are positive constants $K'$ and $\alpha' \leq 1$
such that $f_{\ep }$ is
$C^{1+\alpha' }$ and
the $\alpha'$-H\"older
constant of $f'_{\ep}$ is less than $K'$ for any
$0\leq \ep \leq \ep_{0}$.

\item there are positive constants $K''$ and $\alpha'' \leq 1$
such
that $f_{\ep}'(x)/|x|^{\gamma -1}$ defined on $[-1,0]$
and $f_{\ep}'(x)/|x|^{\gamma -1}$ defined on $[0,1]$
are
$\alpha''$-H\"older continuous
and their
$\alpha''$-H\"older constants
are less than $K''$ for any $0\leq \ep \leq \ep_{0}$,

\item
there are two positive constants $C_{0}$ and $\lambda <1$ such that
the maximum lengths $\lambda_{n, \ep}$ satisfy $\lambda_{n,\ep}\leq
C_{0}\lambda^{n}$ for all positive integers $n$ and all $0\leq \ep \leq
\ep_{0}$.
\end{enumerate}

A function $s$ defined on ${\cal C}^{*}$ is called \underline {H\"older
continuous} if there are two positive constants $C$ and $\lambda <1$
such that
$|s(a^{*})-s(b^{*})|\leq C\lambda^{n}$
for any
$a^{*}$ and
$b^{*}$
in ${\cal C}^{*}$ with
the same first $n$ coordinates.

\vskip5pt
Let ${\cal
A}$ stand for the countable set of points in ${\cal C}^{*}$ whose
coordinates are eventually all
zeroes and let ${\cal B}$ stand for the
complement of ${\cal A}$ in ${\cal C}^{*}$.

\vskip5pt
We can now state the main results of this paper more
precisely.

\vskip5pt
\noindent {\sc Theorem A.} {\em Suppose
$\{ f_{\varepsilon } \}_{0\leq \varepsilon \leq \ep_{0}}$ is a good
family.  There is a family of H\"older
continuous functions $\{ s_{\varepsilon }
\}_{0< \varepsilon \leq \ep_{0}}$
on the dual Cantor set ${\cal C}^{*}$
such that $s_{\varepsilon }$ is the
scaling function of $f_{\varepsilon }$
and

\vskip5pt
\noindent {\sl (1)} for every $0<\varepsilon_{1}\leq \ep_{0}$,
$s_{\varepsilon
}$ converges to $s_{\varepsilon_{1}}$ uniformly on ${\cal
C}^{*}$
as $\varepsilon $ tends to
$ \varepsilon_{1}$,

\vskip8pt
\noindent {\sl (2)} for every $a^{*}\in {\cal C}^{*}$, the limit
$s_{0}(a^{*})$ of
$\{ s_{\varepsilon}(a^{*})\}_{0<\varepsilon \leq \ep_{0}}$ exists as
$\varepsilon $ tends to $0$, the limiting function $s_{0}(a^{*})$
is the scaling function of $f_{0}$ and satisfies:

\vskip5pt
\noindent {\sl (2.1)} $s_{0}$ has jump discontinuities at all points in
${\cal A}$,

\vskip5pt
\noindent {\sl (2.2)} $s_{0}$ is continuous at all points in
${\cal B}$ and the restriction of $s_{0}$ to ${\cal B}$ is a H\"older
continuous function.}

\vskip5pt
Our proof of this theorem depends on a distortion lemma (Lemma 12).
We call it
the \underline {uniform
$C^{1+\alpha}$-Denjoy-Koebe distortion lemma}.

\vskip5pt
More generally, every $f$ on the boundary of
hyperbolicity ${\cal BH}$ has a scaling function (see Figure 6). 
In fact, we show the following theorem.

\begin{figure}[htp]
	\centerline{\psfig{figure=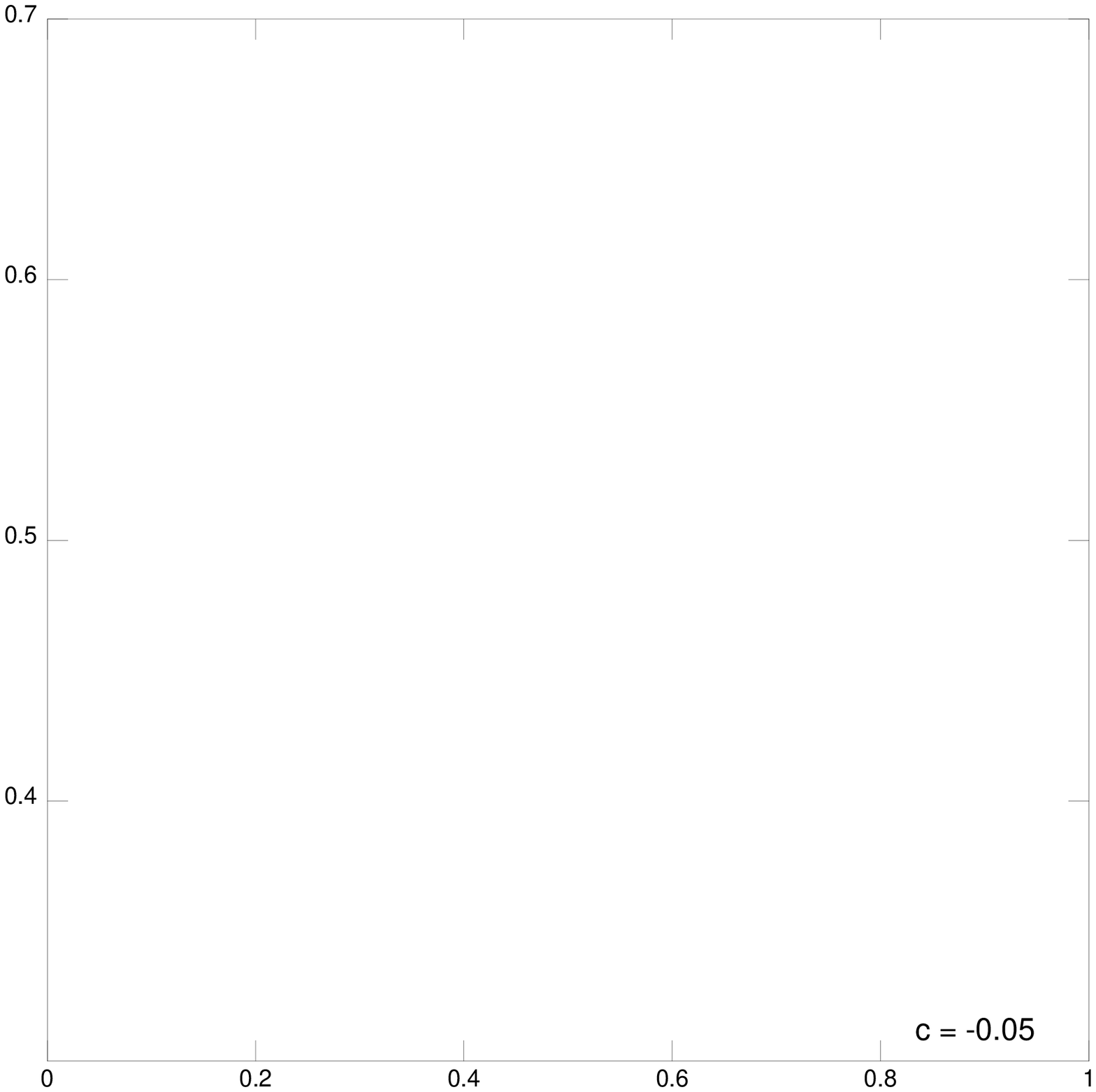,width=.49\textwidth}
		\hfil
		    \psfig{figure=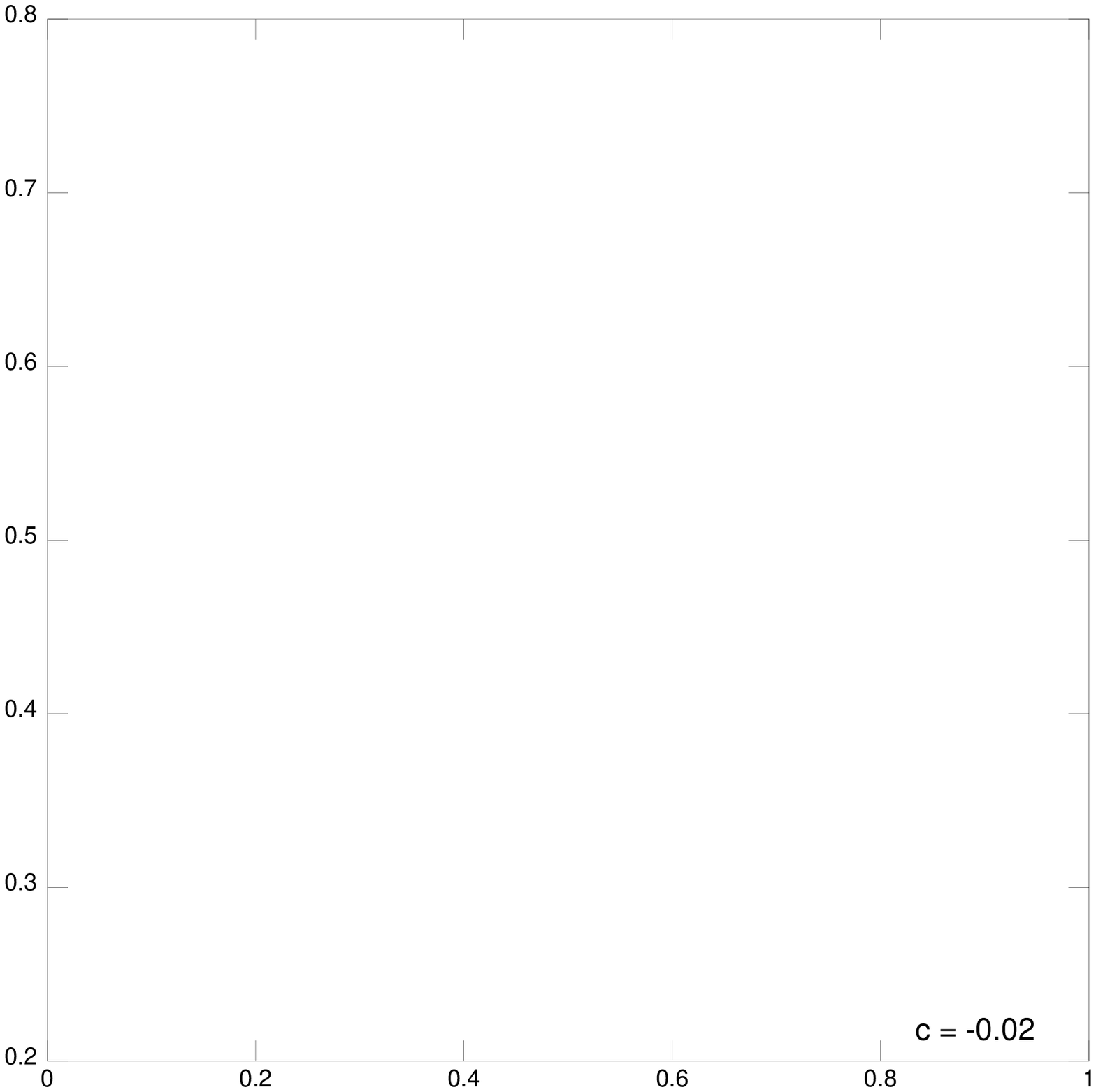,width=.49\textwidth}}
        \centerline{\psfig{figure=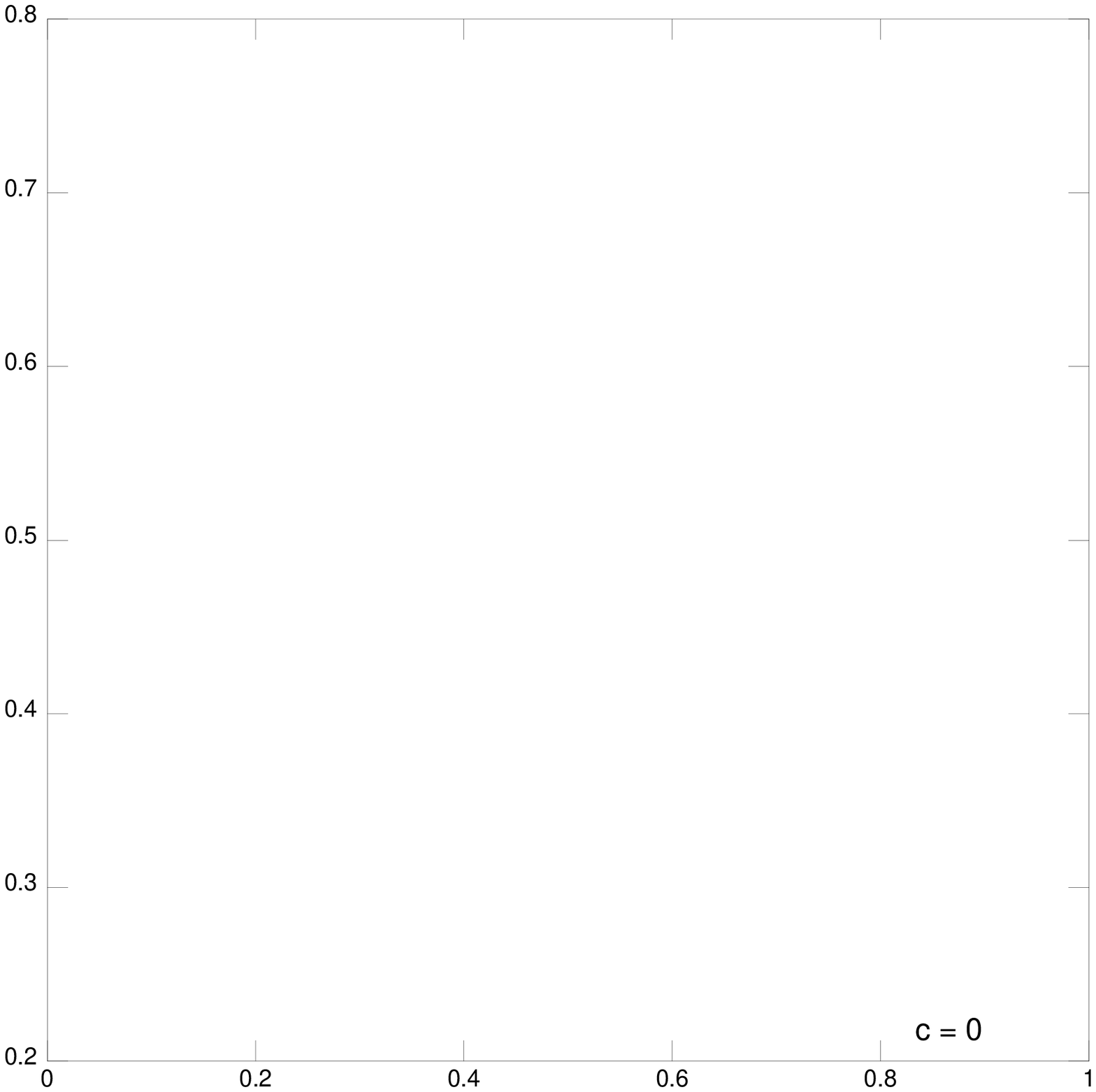,width=.49\textwidth}
		\hfil
		    \psfig{figure=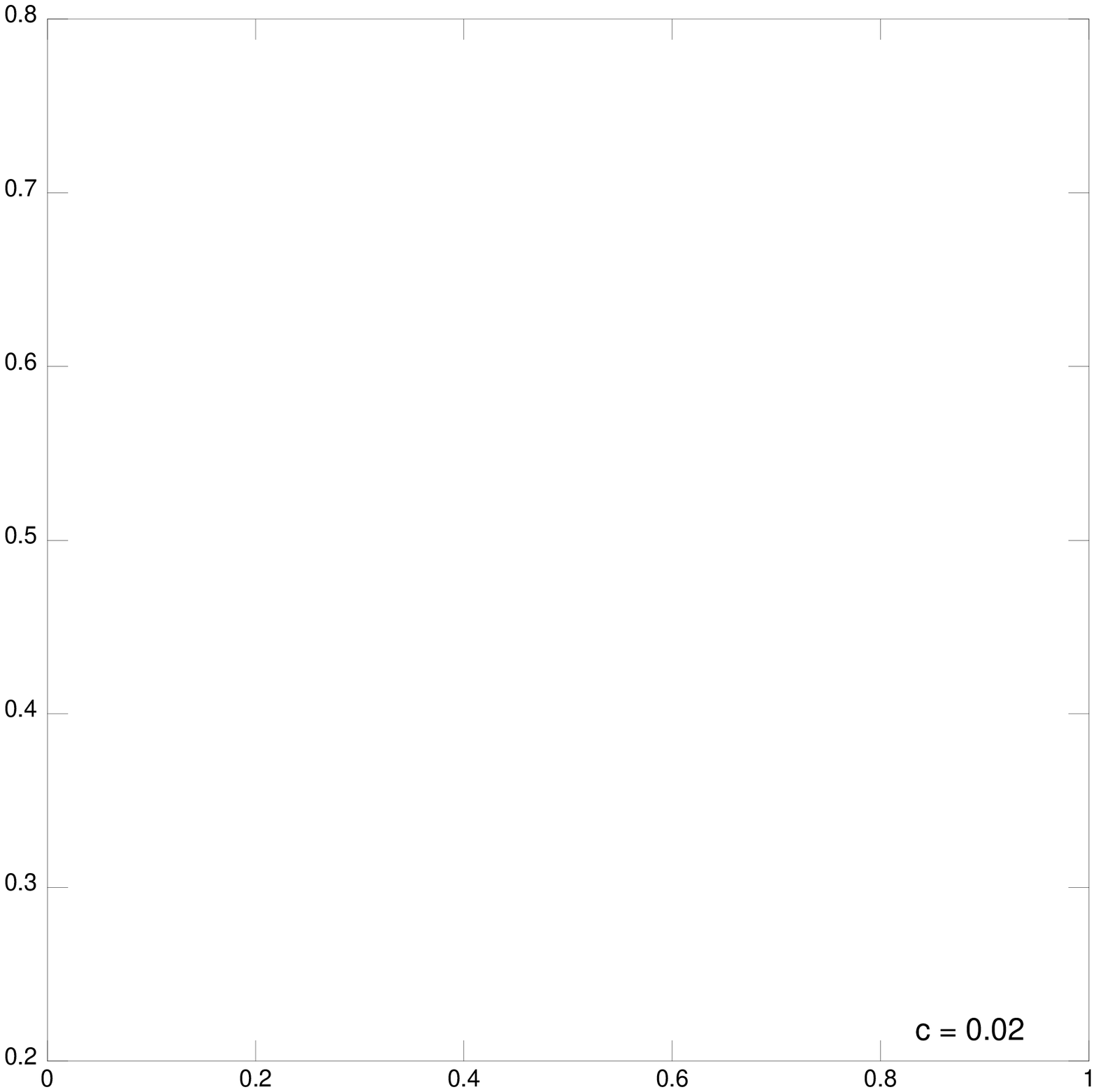,width=.49\textwidth}}
\medskip	
	 \centerline{\sl The graphs of the scaling functions for} 
	 \centerline{\sl $f_c(x) = -x^2 + 2+ cx^2(4-x^2)$,
		with $c=-0.05$, $-0.02$, $0$, and $0.02$.}
\medskip 
         \centerline {\Large \bf Figure 6}
\end{figure}

\vskip5pt
\noindent {\sc Theorem C.} {\em Suppose $f$ is on ${\cal BH}$ and
$\tilde{f}$ is again $f$ viewed in the singular metric associated to
$f$.  There exist the scaling function $s_{f}$ of $f$ and the scaling
function $s_{\tilde{f}}$ of $\tilde{f}$ and these scaling functions
satisfy:

\vskip5pt
\noindent (a) $s_{\tilde{f}}$ is H\"older continuous on ${\cal C}^{*}$,

\vskip5pt
\noindent (b)
$s_{f}$ has jump discontinuities at all points in ${\cal A}$ and
$s_{f}$ is continuous at all points in ${\cal B}$,

\vskip5pt
\noindent (c) the restriction of $s_{f}$ to ${\cal B}$ equals the
restriction of $s_{\tilde{f}}$ to ${\cal B}$.}

\vskip5pt
One example
of a result about a scaling function for a mapping on the boundary of
hyperbolicity is given in the following proposition (see Figure 6, c=0).

\vskip5pt
\noindent {\sc Proposition 2.} {\em
Let $q(x)=1-2x^{2}$ be the mapping of the interval $[-1,1]$.  Then
$s_{q}(a^{*})=1/2$ for all $a^{*}$ in ${\cal B}$ and
$s_{q}(a^{*})\neq 1/2$ for all $a^{*}$ in ${\cal A}$.}

\vskip5pt
Suppose $f_{\ep}$
is a mapping whose graph looks like the graph shown in Figure 2 and
$\{ \eta_{n, \ep }\}_{n=0}^{\infty}$ is the sequence determined by
$f_{\ep}$.  We suppress the subscript $\ep$ when there can be no
confusion.
For every positive integer $n$ and  $I_{w}$ in
$\eta_{n}$, let $ I_{w0}$ and
$I_{w1}$
be the two intervals in
$\eta_{n+1}$ which are contained in
$I_{w}$.
We call the complement of
$I_{w0}$ and
$I_{w1}$ in $I_{w}$ the \underline {gap} on $I_{w}$ and
denote it by
$G_{w}$.
Let $G$ be the complement of $I_{0}$ and
$I_{
1}$. We call $G$ the \underline {leading
gap} and
the set of ratios $\{ |G_{
w}|/|I_{w}|\}$ for all finite strings $w$
of zeroes and ones
the \underline {gap geometry} of the maximal invariant set of
$f_{\ep}$.
We
study the asymptotic dependence on $\ep$ of the gap geometry of the
family of the maximal invariant sets of $f_{\ep}$ for $0\leq \ep \leq
\ep_{0}$.

\vskip5pt
Suppose $\beta $ is a function
defined on
$[0, 1]$.  We say $\beta$ determines asymptotically the gap geometry of
the maximal invariant sets of $f_{\ep}$ for $0\leq \ep \leq \ep_{0}$,
if there is a positive constant $C$ such that for all
$0\leq \ep \leq \ep_{0}$ and all finite strings $w$ of zeroes and ones,

\vskip5pt
\noindent {\sl (1)} $C^{-1}\beta (\ep )\leq
|G_{\ep , w}|/|I_{\ep ,w}|	 \leq C\beta (\ep )$ and

\vskip5pt
\noindent {\sl (2)} $|I_{\ep , wi}|/|I_{\ep ,w}|\geq C^{-1}$, where
$i$ is either one or zero.

\vskip5pt
\noindent The constant $C$ is called a determining constant.

\vskip5pt
Suppose
$\{ f_{\varepsilon } \}_{0\leq \ep \leq \ep_{0}}$ is a good family.
Then the sizes of the leading gaps $G_{\ep}$
is of order of
$\ep^{\frac{1}{\gamma}}$.  Moreover, we
prove the following theorem.

\vskip5pt
\noindent {\sc Theorem B.} {\em The family of the maximal invariant
sets of $f_{\ep}$ for $0<\ep \leq \ep_{0}$ is a family of
Cantor sets. Furthermore,
the function $\ep^{\frac{1}{\gamma}}$ determines asymptotically the gap
geometry of the family of maximal invariant sets of $f_{\ep}$ for
$0< \ep \leq \ep_{0}$.}

\vskip5pt
Suppose $HD(\ep)$ is the Hausdorff dimension of the maximal
invariant set of $f_{\ep}$ and
$s_{\ep}$ is the scaling
function of $f_{\ep}$ for $0\leq \ep \leq \ep_{0}$.
Theorem B has the following two corollaries.

\vskip5pt
\noindent {\sc Corollary 2.}
{\em There is a
positive constant $C$ which does not depend on $\ep$ such that
\[0<  HD(\ep) \leq 1-C\ep^{\frac{1}{\gamma}}\]
for all $0\leq \ep \leq \ep_{0}$.}

\vskip5pt
\noindent {\sc Corollary 3.}
{\em
There is a positive
constant $C$ which does not depend on $\ep$ such that
\[ 1-C^{-1} \varepsilon^{\frac{1}{\gamma}} <
s_{\varepsilon}((a^{*}0.))
+s_{\varepsilon}((a^{*}1.)) <1-C \varepsilon^{\frac{1}{\gamma }}\]
for all $a^{*}\in {\cal C}^{*}$ and all $0\leq \ep \leq \ep_{0}$.}

\vskip5pt
For the family $\{ 1+\ep -(2+\ep
)x^{2}\}_{0\leq \ep \leq 1}$ in ${\cal H}\cup {\cal BH}$, we can say
even more:

\vskip5pt
\noindent {\sc Proposition 3.}
{\em
There is a constant $C>0$ which does not depend on
$\varepsilon $ such that
\[ 1-C^{-1}\sqrt{\varepsilon } \leq HD(\varepsilon ) \leq
1-C\sqrt{\varepsilon }\]
for all $0\leq \varepsilon \leq 1$.}

\pagebreak
\noindent {\Large \bf \S 2 The Boundary Of Hyperbolicity, ${\cal BH}$}

\vskip5pt
The example studied by Ulam and von Neumann [UN] in 1947 is
the \nl\ \tr\ $q(x)=1-2x^{2}$. They discovered the
density \fn\ $\rho_{q}(x)=
1/(\pi \sqrt {1-x^{2}})$ of the unique \acim\ for this mapping. From the metric point of view, $dy =2\rho_{q}(x)dx $ is a
singular
metric on \intl .  Under the corresponding change of coordinate
$y=h(x)$,
$q$ becomes a piecewise linear \tr\ $\tilde{q}(y)=h\circ q\circ h^{-1}(y)=
1-2|y|$. The point is that the dynamics of $\tilde{q}$ is more easily
understood.

\vskip10pt
\noindent {\bf \S 2.1	The singular change of metric on the
interval.}

\vskip5pt
Suppose $f$ is a $C^{1}$ self mapping of an interval
with a unique \cp\ $c$.
We always
make the following assumptions: {\sl (1)} $f$ is a $C^{1}$ self mapping
of $[-1,1]$,  {\sl (2)} $f$
is increasing
on $[-1,c]$ and decreasing on $[c,1]$,	{\sl (3)} $f$ maps $c$ to $1$ and
{\sl (4)} $f$ maps $-1$ and $1$ to $-1$ (see Figure 1a).   Without loss
of generality, we always assume $c$ equals $0$.

\vskip5pt
Let
$r_{f}(x) = f'(x)/|x|^{\gamma -1}$ for $x\neq 0$.
We
say
$f$ has \pl\ at the \cp\ if there is some number
$\gamma >1$ such that the limits of $
r_{f}(x) $ as $x$ increases to zero and
as $x$ decreases to zero
exist and equal nonzero numbers $A$ and $-B$, respectively.
For
example, $f$ has power law at the critical point if
$f(x)=1-(A/\gamma)|x|^{\gamma }$ for negative $x$ close to zero
and $f(x)= 1-(B/\gamma)|x|^{\gamma }$ for positive $x$ close to zero.
We call the ratio of
$A$ to $-B$, which is the limit of $f'(-x)/f'(x)$ as $x$ decreases to
zero,
the \underline {asymmetry} 
of $f$ at the critical point (see [J3]).
We always assume that $f$ has power law $|x|^{\gamma} $ for some
$\gamma >1$ at the critical point.

\vskip7pt
We define the \underline {singular metric associated to
$f$} to be
\[ dy =\frac{dx}{(1-x^{2})^{\frac{\gamma-1}{\gamma}}}\]
on $[-1,1]$
The corresponding change of coordinate on $[-1,1]$ is $y=h_{\gamma}(x)$,
where
\[ h_{\gamma }(x) =-1+b
\int_{-1}^{x}
\frac{dx}{(1-x^{2})^{\frac{\gamma-1}{\gamma}}}\]
with $b=2/\int_{-1}^{1} \sms\ $.
The representation of $f$ under the singular metric associated to
$f$ is
\[ \tilde{f}= h_{\gamma}\circ f\circ h_{\gamma}^{-1}.\]

\vskip7pt
Before we see some properties of $f$ and $\tilde{f}$, we state some
definitions.

\vskip5pt
\noindent {\sc Definition 1.} {\em
Suppose $I$ is an interval and $g$
is a function on $I$. We say that

\vskip5pt
\noindent (a) $g$ is an \underline {embedding} if $g$ is
a homeomorphism from $I$ to $g(I)$,

\vskip5pt
\noindent (b) $g$ is a
\underline {$C^{1}$-embedding} if $g$ has a
continuous derivative $g'$ on $I$ and the derivative of
$g$ at
every point in $I$ is not zero,  as usual, the derivative of $g$
at each boundary point of $I$ is a one-sided limit,

\vskip5pt
\noindent (c) $g$ is a \underline {$C^{1+\alpha}$-embedding} for some
$0< \alpha \leq 1$ if $g$ is a $C^{1}$-embedding
and the derivative $g'$ on $I$ is $\alpha$-H\"older continuous on $I$.

\vskip5pt
\noindent If $\alpha =1$, we usually say $g$ is a
$C^{1,1}$-embedding in (c).}

\vskip7pt
Suppose $I$ is an interval and $g$
is an $\alpha$-H\"older continuous function
on $I$ for some
$0<\alpha \leq 1$.
There is
a positive constant $K$ such that $|g(x)-g(y)|\leq K|x-y|^{\alpha }$ for
all $x$ and $y$ in $I$.  The smallest such $K$ is called the
$\alpha$-H\"older
constant of $g$.  If $\alpha =1$, the smallest such $K$ is usually
called the Lipschitz constant of $g$.

\vskip7pt
\noindent {\sc Lemma 1.} {\em The mapping $\tilde{f}$
is continuous on $[-1,1]$ and the
restrictions of $\tilde{f}$ to
$[-1,0]$ and to $[0,1]$ are $C^{1}$-embeddings.}

\vskip7pt
\noindent {\sl Proof.} If $y$ is not one of $0$, $1$ and
$-1$, then $\tilde{f}$
is differentiable at $y$. Suppose $x$ is the
preimage of $y$ under $h_{\gamma}$.  By the chain rule,
\[ \tilde{f}^{\prime }(y) = f'(x)
(1-x^{2})^{\frac{\gamma-1}{\gamma}}/
(1-(f(x))^{2})^{\frac{\gamma-1}{\gamma}}. \hskip100pt (EQ \hskip5pt
1.1).\] Using this equation, we can get that
$\tilde {f}'(0-)$ and
$\tilde {f}'(0+)$ exist and equal nonzero numbers and that
$\tilde {f}'(-1)=
(f'(-1))^{\frac{1}{\gamma}}$ and
$\tilde {f}'(1)=-
(|f'(1)|)^{\frac{1}{\gamma}}$.
\Q

\vskip7pt
\nd {\sc Remark 1.}
The inverse
of $h_{\gamma}$
is $C^{1}$.
If the restrictions of
$r_{f}$ to $[-1, 0)$ and to $(0, 1]$ are $\alpha$-H\"older
continuous for some $0< \alpha \leq 1$,
then the restrictions of
$\tilde{f}$ to $[-1,0]$ and to $[0,1]$ are at least $C^{1+\alpha}$
because of (EQ 1.1) (see [J4]).

\vskip7pt
\noindent{\sc Lemma 2.} {\sl Suppose $\tilde{f}$ is a continuous
self mapping of $[-1, 1]$.	Assume $0$ is
the unique turning point, $f$ maps $0$ to $1$ and maps $-1$ and $1$ to
$-1$ and the restrictions
of $\tilde{f}$ to $[-1,0]$ and to $[0,1]$ are $C^{1}$-embeddings.	Then
$f=h_{\gamma}^{-1}\circ \tilde{f}\circ h_{\gamma}$ is a $C^{1}$ mapping
and
has the power law $|x|^{\gamma }$ at the critical
point $0$ for any $\gamma >1$.}

\vskip7pt
\noindent{\sl Proof.}
If $x$ is not one of $0$, $1$ and $-1$, then
$f$
is differentiable at $x$.  Suppose $y=
h_{\gamma}(x)$.
By the chain rule,
\[ f^{\prime }(x) = \tilde{f}'(y)
(1-(h_{\gamma}^{-1}\circ \tilde{f}(y))^{2})^{\frac{\gamma-1}{\gamma}}
/(1-(h_{\gamma}^{-1}(y))^{2})^{\frac{\gamma-1}{\gamma}}.
\hskip50pt (EQ \hskip5pt 1.2).\]
Using this equation,
$f'(-1)=
(\tilde{f}'(-1))^{\gamma}$ and
$f'(1)=
-|\tilde{f}'(1)|^{\gamma}$,  and the limits of $
r_{f}(x)$ as $x$ increases to zero and as
$x$ decreases to zero exit and
equal nonzero numbers.
\Q

\vskip7pt
\nd {\sc Remark 2.} The mapping $h_{\gamma}$ is
$(1/\gamma)$-H\"older
continuous.
If
the restrictions of $\tilde{f}$ to $[-1,0]$ and
to $[0,1]$ are
$C^{1+\alpha }$ embeddings for some $0< \alpha \leq 1$, then
$f$ is $C^{1+\frac{\alpha }{\gamma}}$ and the
restrictions of $r_{f}$ to $[-1,
0)$ and to $(0,1]$ are $\alpha /\gamma$-H\"older continuous
because of (EQ 1.2) (see [J4]).

\vskip10pt
\noindent {\bf \S 2.2 The definition of the boundary of hyperbolicity}

\vskip5pt
Let $f_{0}$ and $f_{1}$ be the restrictions of $f$ to
$[-1,0]$ and to $[0, 1]$, respectively. Then $f_{0}$
and $f_{1}$ are both embeddings.
Let $g_{0}$ and $g_{1}$ be the inverse of $f_{0}$ and
$f_{1}$. For a finite string $w=i_{n}\cdots i_{0}$ of zeroes and ones,
let $g_{w}$
be the composition, $g_{w}=
g_{i_{0}}\circ \cdots \circ g_{i_{n}}$.  Define $I_{w}$ to be the image
of $[-1,1]$ under $g_{w}$.
The \npt\	of $[-1, 1]$ determined by $f$ is the
collection
of $I_{w}$ for all finite strings $w$ of zeroes and ones
of
length $n+1$. We denote it by $\eta_{n,f}$ or just by $\eta_{n}$ when
there is no possibility for confusion.
We use $\lambda_{n}$ to denote the maximum length of the intervals in
$\eta_{n}$.  The $n^{th}$-partition $\eta_{n,\tilde{f}}$ of $[-1,1]$
determined by $\tilde{f}$ and the maximum length
$\lambda_{n,\tilde{f}}$ of the intervals in $\eta_{n,\tilde{f}}$ are
defined similarly.

\vskip5pt
\noindent {\sc Definition 2.} {\em
We say that
the sequence of nested partitions $\{ \eta_{n}\}_{n=0}^{\infty}$
determined by $f$ decreases
exponentially if there are two
positive constants $C$ and $\lambda <1$ such that $\lambda_{n} \leq
C\lambda^{n}$ for all positive integers $n$.}

\vskip5pt
\noindent{\sc Lemma 3.} {\em
The sequence of nested partitions
determined by $f$ decreases
exponentially if and only if
the sequence of nested partitions
determined by $\tilde{f}$ decreases
exponentially.}

\vskip5pt
\noindent {\sl Proof:} Because $h_{\gamma}$ is $(1/\gamma)$-H\"older
continuous and
the inverse of $h_{\gamma}$ is $C^{1}$, we can easily see this
lemma.
\Q

\vskip5pt
\noindent{\sc Definition 3.} {\em The \nl\ mapping $f$ is on the
\underline {boundary of hyperbolicity},  ${\cal BH}$, if

\vskip5pt
\nd (a) the restrictions of $\tilde{f}$ to $[-1,0]$ and to $[0,1]$ are
$C^{1+\alpha }$ embeddings for some $0< \alpha \leq 1$,

\vskip5pt
\nd (b) the sequence of nested partitions $\{ \eta_{n} \}$ determined by
$f$ decreases exponentially. }

\vskip5pt
Next lemma follows from Remark 1 and Remark 2.

\vskip5pt
\nd {\sc Lemma 4.} {\em
The \nl\ mapping $f$ is on ${\cal BH}$ if and only if

\vskip5pt
\nd $(i)$ $f$ is
$C^{1+\alpha }$ for some $0< \alpha \leq 1$ and
the restrictions of $r_{f}$ to $[-1, 0)$ and to $(0, 1]$ are
$\beta$-H\"older continuous
for some $0< \beta \leq 1$ and

\vskip5pt
\nd $(ii)$
the sequence of nested partitions determined by $f$ decreases
exponentially. }

\vskip5pt
The condition $(i)$ in Lemma 4 and
$r_{f}(0-)=r_{f}(0+)$ are equivalent to the statement that $f(x)
=F(-|x|^{\gamma})$
where $F$ is a $C^{1+\alpha }$ diffeomorphism from $[-1,0]$ to $[-1,1]$
(see [J4] for more details).

\vskip5pt
We give two examples of mappings on ${\cal BH}$.
The reader may refer to the paper [J2]
for the proofs that these two examples are on ${\cal BH}$

\vskip5pt
\nd {\sc Example 1.} {\em  Mappings $f$ such that $(1)$ $f$ is
$C^{3}$ with nonpositive
\sd , $(2)$ $f$ is expanding at both boundary points of
$[-1,1]$, that is, $f'(-1)$ and $|f'(1)|$ are greater than one, and
$(3)$ the restrictions of $r_{f}$ to $[-1,0)$ and to $(0, 1]$ are
$\alpha$-H\"older continuous
for some $0< \alpha \leq 1$.}

\vskip5pt
The \sd\ $S(f)$
of $f$ is
$S(f) =f'''/f'-(3/2) (f''/f')^{2}$.

\vskip5pt
\nd {\sc Example 2.} {\em Mappings $f$ such that
$(1)$ the restrictions of $\tilde{f}$ to $[-1,0]$ and
to $[0,1]$ are $C^{1,1}$ embeddings and
$(2)$ all the periodic points of $f$ are expanding,
that is, the absolute values of the eigenvalues of $f$ at all
periodic points are greater than one.}

\vskip5pt
The eigenvalue of $f$ at a periodic point
$p$ of period $n$ of $f$ is $e_{f}(p)= (f^{\circ
n})'(p)$.

\vskip20pt
\noindent {\Large \bf \S 3 The Space Of Hyperbolic Mappings, ${\cal H}$}

\vskip5pt
A type of perturbation of
$q:x \mapsto 1-2x^{2}$
is a mapping $q_{\varepsilon}:x \mapsto 1+\varepsilon -(2+\varepsilon
)x^{2}$ for a \pn\ $\varepsilon$. The mapping
$q_{\varepsilon}$ maps the critical point out of $[-1,1]$ and
it does not keep $[-1,1]$ invariant but invariant a Cantor set which has
bounded geometry.
In this paper, we study the asymptotic
behavior of certain mappings like
$q_{\varepsilon}$ as they approach the boundary of hyperbolicity (see
Figure 5).

\vskip10pt
\nd {\bf \S 3.1 The definition of the scaling function.}

\vskip5pt
Suppose $\ep $ is a positive number and $f_{\varepsilon}$ is a $C^{1}$
mapping from $[-1,1]$ to the real line with a unique \cp\ $c$.
We always make the following assumptions:
(1) $f_{\varepsilon}$ is increasing
on $[-1,c]$ and decreasing on $[c,1]$, (2) $f_{\ep}(c) =1+\ep$ and (3)
$f_{\ep}$ maps $1$ and $-1$ to $-1$
(see Figure 2).  Without loss of generality, we always assume $c$
equals $0$.

\vskip5pt
Let $f_{\ep , 0}$ and $f_{\ep , 1}$ be the restrictions of $f_{\ep}$
to
$[-1,0]$ and to $[0, 1]$. They are two
embeddings.
Let $g_{\ep ,0}$ and $g_{\ep ,1}$ be the inverses of
$f_{\ep, 0}$ and $f_{\ep, 1}$.	For a finite string $w=i_{0}\cdots i_{n}$
of zeroes and ones,  let $g_{\ep , w}$ be
the composition
$g_{\ep , w}= g_{\ep ,i_{0}}\circ \cdots \circ g_{\ep ,i_{n}}$ and
$I_{\ep ,w}$ be
the image of $[-1, 1]$ under $g_{\ep, w}$.
Suppose $\eta_{n,\ep}$ is
the collection
of $I_{\ep ,w}$ for all finite strings $w$ of zeroes and ones
of length $n+1$ and $\lambda_{n,\ep }$ is the maximum
length of the intervals in
$\eta_{n,\ep}$.  The union of the intervals in $\eta_{n,\ep}$ covers the
maximal invariant set of $f_{\ep}$.
We always assume that $\lambda_{n, \ep}$ goes to zero as n increases to
infinity.

\vskip5pt
For every interval in $\eta_{n, \ep}$, there is the labelling
$w$ where $w$
is the finite string of zeroes and ones such that this interval is the
image $I_{w}$ of $[-1, 1]$ under $g_{\ep , w}$.
There are two topologies
on the set of all the labellings $w$ .	One topology is induced by
reading the labellings $w$ from left to right; the other topology is
induced by reading the labellings $w$ from right to left.

\vskip5pt
Suppose we read all the labellings $w$ from left to right and ${\cal
C}_{n} = \{ w_{n} |$ $ w_{n}=(.i_{0}i_{1}\cdots i_{n}),$ where $i_{k}$
is either $0$ or $1$ for $k\geq 0$ and $n \geq 0\} $.  Let ${\cal
C}_{n}$ have the product topology.  The continuous mapping $\sigma_{n}:
{\cal C}_{n+1} \mapsto {\cal C}_{n}$ is defined by
$\sigma_{n}((.i_{0}i_{1}\cdots i_{n})) =(.i_{1}\cdots i_{n})$ for $n\geq
0$.
The
pairs $\{ ({\cal C}_{n}, \sigma_{n}) \}_{n =0}^{\infty} $ form
an inverse limit set.	Let
${\cal C}$ be the inverse limit of this inverse limit set and
$\sigma $ be the induced mapping on ${\cal C}$.  We call ${\cal C}$ the
\underline {topological Cantor set}.  For any $a$ in ${\cal C}$, it is
an infinite string
of zeroes and ones extending to the right, that is,
$a =(.i_{0}i_{1}\cdots
)$ where $i_{k}$ is either zero or one for $k\geq 0$.  The
mapping $\sigma$ is the shift
mapping on ${\cal C}$, that is, $\sigma $ maps $(.i_{0}i_{1}\cdots
)$ to $(.i_{1}\cdots )$.  We call $({\cal C}, \sigma )$ the
symbolic dynamical system of $f_{\ep}$ because of the following lemma.

\vskip5pt
\nd {\sc Lemma 5.} {\em
Suppose $\Lambda_{\ep}$
is the maximum invariant set of $f_{\ep}$.
There is a homeomorphism $h_{\ep}$ from ${\cal C}$ to $\Lambda_{\ep}$
such
that $h_{\ep}\circ \sigma  = f_{\ep}\circ h_{\ep}$. In the
other words,
$(\Lambda_{\ep}, f_{\ep})$ and $({\cal C},$ $\sigma )$ are conjugate.}

\vskip5pt
\nd {\sl Proof.} Suppose $a=(.i_{0}i_{1}\cdots )$ is any point in ${\cal
C}$.   Let $w_{n}=(.i_{0}\cdots i_{n})$ be the first $n+1$
coordinates of $a$.
The intersection of nested intervals $\{ I_{w_{n}}\}_{n=0}^{\infty} $
is nonempty and contains only one point $x(a)$ because the length of
$I_{w_{n}}$
goes to zero as $n$ increases to infinity.
Define $h_{\ep }(a) =x(a)$. Then $h_{\ep}$ is a
homeomorphism from ${\cal C}$ to $\Lambda_{\ep}$ and
$h_{\ep}\circ \sigma  = f_{\ep}\circ h_{\ep}$. \Q

\vskip5pt
Suppose we read all the labellings $w$ from right to left and ${\cal
C}^{*}_{n} = \{ w^{*}_{n} |$ $ w^{*}_{n}=(i_{n}\cdots i_{1}i_{0}.)$,
where $i_{k}$ is
either zero or one and $n \geq 0\}$.  Let ${\cal
C}^{*}_{n}$ have the product topology.	The continuous mapping
$\sigma^{*}_{n}: {\cal C}^{*}_{n+1} \mapsto {\cal C}^{*}_{n}$ is defined
by
$\sigma^{*}_{n}((i_{n}\cdots i_{1}i_{0}.)) =(i_{n}\cdots i_{1}.)$ for
$n\geq 0$.
The
pairs $\{ ({\cal C}^{*}_{n}, \sigma^{*}_{n}) \}_{n =0}^{\infty}
$ also form an inverse limit set.   Let
${\cal C}^{*}$ be the inverse limit of this inverse limit set and
$\sigma^{*} $ be the induced mapping on ${\cal C}^{*}$.
We call ${\cal C}^{*}$ the
\underline {dual Cantor set}.  Any $a^{*}$ in ${\cal C}^{*}$ is an
infinite string of zeroes and ones extending to the left, that
is,	$a^{*}=(\cdots
i_{1}i_{0}.)$ where $i_{k}$ is either zero or
one for $k\geq 0$. 	We call $a^{*}$ a ``dual point'' of $f_{\ep}$.
The mapping $\sigma^{*}$ is the shift mapping on ${\cal C}^{*}$, that
is, $\sigma^{*} $ maps $(\cdots i_{1}i_{0}.)$ to
$(\cdots i_{1}.)$.  We call $({\cal C}^{*},
\sigma^{*})$ the dual
symbolic dynamical system of $f_{\ep}$.

\vskip5pt
A sequence $\{ x_{n}\}_{n=0}^{\infty }$ in
the maximum invariant set of $f_{\ep}$ is a sequence of
backward images of $x_{0}$ under $f_{\ep }$ if
$f_{\ep}(x_{n})=x_{n-1}$ for all positive integers $n$.  The dual Cantor
set will not represent the maximal invariant set of
$f_{\ep}$, but there is a one-to-one
corresponding from the dual Cantor set to the set of sequences
of backward images of $x_{0}$ under $f_{\ep}$ for all points
$x_{0}$ in the maximal invariant set of $f_{\ep}$.
The
\underline {scaling function} of $f_{\ep}$ is defined on the dual Cantor
set ${\cal C}^{*}$ if it exists (see [S2] and [J3]). For the sake of
completeness, we give the definition of scaling function as the following.

\vskip5pt
Suppose $a^{*}$ is in ${\cal C}^{*}$, so that $a^{*}$ is an infinite
string
of zeroes and ones extending to the left.  Assume $a^{*}= (\cdots
wi.)$
where $w$ is a finite string of zeroes and ones and $i$ is either zero
or one.  Note that $I_{wi}$ is a
subinterval of $I_{w}$.  Let $s(wi)$ equal the ratio of the lengths,
$|I_{wi}|/ |I_{w}|$.  We let $s(a^{*})$ be the limit set of $s(wi)$ as
the length of $w$ tends to infinity. 

\vskip5pt
\nd {\sc Definition 4.} {\em Suppose $f_{\ep}$ is a mappings in ${\cal
H}$.  If the limit set $s(a^{*})$ consists of only one number for 
$a^{*}$ in ${\cal C}^{*}$, then we say there is the scale $s(a^{*})$ of
$f_{\ep}$ at $a^{*}$. If there is the scale $s(a^{*})$ of $f_{\ep}$ for
every $a^{*}$ in $C^{*}$, then we call $s:C^{*}\mapsto {\bf R}^{1}$ 
the scaling function of $f_{\ep}$. Note that the scaling function 
$s(a^{*})$ of
$f_{\ep}$ depending
on $\ep$. Sometimes we denote it by $s_{\ep}(a^{*})$.}

\vskip5pt
\nd {\sc Remark 3.} For $f$ on ${\cal BH}$, we can use the
the same arguments as Definition 4 to
define the scaling function $s_{f}$ of $f$ on
${\cal C}^{*}$ by the
sequence of nested partitions $\{ \eta_{n}\}_{n=0}^{\infty}$ determined
by $f$.

\vskip10pt
\nd {\bf \S 3.2 The definition of the space of hyperbolic mappings}

\vskip5pt
Suppose
$\{ \eta_{\ep , n}\}_{n=0}^{\infty}$
is the sequence determined by $f_{\ep}$.
Just as in Definition 2, we say
the sequence
$\{ \eta_{n,\ep}\}_{n=0}^{\infty}$ determined by $f_{\ep}$ decreases
exponentially if $\lambda_{n, \ep}$ decreases exponentially.

\vskip5pt
\noindent{\sc Definition 5.} The \nl\ mapping $f_{\ep}$
is in the space of hyperbolic mappings, ${\cal
H}$, if

\vskip5pt
\nd (a) $f_{\ep}$
is
$C^{1+\alpha }$  for some $0< \alpha \leq 1$ and

\vskip5pt
\nd (b)
the sequence $\{ \eta_{n, \ep }\}_{n=0}^{\infty }$ determined by
$f_{\ep}$ decreases exponentially.

\vskip5pt
We give two examples of mappings in ${\cal H}$.  They are similar to
Example 1 and Example 2  in \S 2.

\vskip5pt
\nd {\sc Example 3.} {\em Mapping $f_{\ep}$ such that
$(1)$ $f_{\ep}$ is a $C^{3}$ mapping on $[-1, 1]$ with
nonpositive \sd\ and $(2)$ $f$ is expanding at both boundary points of
$[-1,1]$, that is,
$f_{\ep}'(-1)$ and $|f_{\ep}'(1)|$ are greater than
one.}

\vskip5pt
\nd {\sc Example 4.} {\em Mapping $f_{\ep}$ such that $(1)$
$f_{\ep
}$ is $C^{1,1}$ and $(2)$ all the periodic points of $f_{\ep}$ are
expanding, that is, the absolute values of all eigenvalues of
$f_{\ep}$ at periodic points
are greater than one.}

\vskip5pt
The proofs, that Example 3 and Example 4 are in ${\cal H}$, are
similar to the proofs of Example 1 and 2 in [J2].  

\vskip5pt
\nd {\sc Definition 6.} {\em A function $s$ defined on dual
Cantor set $C^{*}$ is \underline {H\"older}
\underline {continuous} if there are two positive constants $C$ and
$\lambda <1$ such that
$|s(a^{*})-s(b^{*})|\leq C\lambda^{n}$ for any
$a^{*}$ and $b^{*}$ in $C^{*}$ with the same first
$n$ coordinates.  We call
$C$ a H\"older constant of $s$.}

\vskip5pt
\nd {\sc Lemma 6.} {\em Suppose $f_{\ep}$ is in ${\cal H}$.
There exists a H\"older continuous scaling
function $s_{\ep }$ of $f_{\ep }$.}

\vskip5pt
\nd {\sl Proof.}
Let $d_{\ep}$ be the minimum value of the restriction of
$f_{\ep}$ to the union of $I_{\ep , 0}$ and $I_{\ep , 1}$.
We suppress
$\ep$ if there can be no confusion. Note
that $d_{\ep}$ goes to zero as $\ep$ decreases to zero.
For any $a^{*}$ in ${\cal C}^{*}$, we use $w_{n}i$
to denotes the
first $(n+1)$ coordinates of $a^{*}$ and  $s(
w_{n}i)$ to denote the ratio, $|I_{
w_{n}i}|/|I_{w_{n}}|$.
By (b) of Definition 5, we have two positive constants $C_{0}$
and $\lambda<1$ such that $\lambda_{n, \ep} \leq C_{0}
\lambda^{n}$.  Let $K$ be the H\"older constant of $f_{\ep}'$ on
$[-1,1]$.
Because $s(w_{m}i) = (|(f^{\circ (m-n)})'(x)|/ |(f^{\circ
(m-n)})'(y)|) s(w_{n}i)$ for some $x$ and $y$ in $I_{w_{m}}$,
by the naive distortion lemma,	there is
a constant
$C_{\ep}$ which equals $C_{0}K/(d_{\ep}(1-\lambda^{\alpha }))$ such
that for any $m>n>0$,
\[ |s(w_{m}i)-
s(w_{n}i)
|\leq C_{\ep } |I_{w_{n}}|^{\alpha }.\]
The last inequality implies that the limit of sequence
$\{ s(w_{n}i) \}_{n=0}^{+\infty }$ exists as the length of $w_{n}i$
increases
to infinity.  We denote this limit by $s(a^{*})$ or $s_{\ep}(a^{*})$ if
we need to indicate dependence on $\ep$.
Let $m$ tend to infinity,  then $|s(a^{*})-
s(w_{n}i)
|\leq C_{\ep } |I_{\ep , w_{n}}|^{\alpha }$ for all positive
integers $n$.

\vskip5pt
Suppose $a^{*}$ and $b^{*}$ are in ${\cal C}^{*}$ with the same
first $(n+1)$ coordinates, that is, $a^{*}= (\cdots w_{n}i.)$ and
$b^{*}=( \cdots w_{n}i.)$.  Because
$|s(a^{*})-
s(w_{n}i)
|\leq C_{\ep } |I_{w_{n}}|^{\alpha }$ and
$|s(b^{*})-
s(w_{n}i)
|\leq C_{\ep } |I_{w_{n}}|^{\alpha}$, we have
that $|s(a^{*})- s(b^{*})|\leq
2C_{\ep } |I_{w_{n}}|^{\alpha }\leq 2C_{0}C_{\ep} \lambda^{n}$. In
other words,
$s_{\ep }$ is H\"older continuous on ${\cal C}^{*}$ with a H\"older
constant $2C_{0}C_{\ep}$. \Q

\vskip5pt
In the next chapter, we will study the asymptotic behavior of scaling
function $s_{\ep}$, as well as the geometry of Cantor set $\Lambda_{\ep}$, 
which is the
maximal invariant set of $f_{\ep}$, for $f_{\ep}$ in ${\cal H}$.

\vskip20pt
\nd {\Large \bf \S 4 Asymptotic Geometry Of Cantor Sets}

\vskip5pt
Suppose $f_{\ep}$ is in ${\cal H}$.  Let $r_{\ep}(x)
= f_{\ep}'(x)/|x|^{\gamma -1}$ for nonzero $x$ in $[-1,1]$.
We say that
$f_{\ep}$ has \pl\ at the \cp\ if there is some
$\gamma >1$ such that the limits of $
r_{\ep}(x) $ as $x$ increases to zero and
as $x$ decreases to zero
exist and equal nonzero numbers $A_{\ep}$ and $-B_{\ep}$,
respectively.

\vskip5pt
We define the \underline {smooth metric associated to
$f_{\ep}$} to be
\[ dy
=\frac{dx}{((1+\varepsilon)^{2}-x^{2})^{\frac{\gamma-1}{\gamma}}} \]
on $[-1,1]$.
The corresponding change of coordinate is $y=h_{\gamma , \ep}$ where
\[ h_{\gamma, \ep }(x)
=-1+b_{\ep} \int_{-1}^{x}
\frac{dx}{((1+\varepsilon)^{2}-x^{2})^{\frac{\gamma-1}{\gamma}
}}\]
with $b_{\ep}=2/\int_{-1}^{1} \sms\
$.
The representation of $f_{\ep}$ under the smooth metric associated to
$f_{\ep}$ is
\[ \tilde{f_{\ep}}= h_{\gamma , \ep}\circ
f_{\ep}\circ h_{\gamma , \ep }^{-1}.\]

\vskip5pt
\noindent {\sc Lemma 7.} {\em If $f_{\ep}$ has power law
$|x|^{\gamma }$ with $\gamma >1$, then the mapping $\tilde{f_{\ep}}$
is continuous on $[-1,1]$ and the restrictions of $\tilde{f}_{\ep}$ to
$[-1, 0]$ and to $[0,1]$ are $C^{1}$ embeddings.}

\vskip5pt
\nd {\sc Remark 4.} The mapping $h_{\gamma , \ep}$ is a $C^{\infty }$
diffeomorphism from $[-1, 1]$ to itself.
If
the derivative
$f_{\ep}'$ and the restrictions of
$r_{\ep}$ to $[-1, 0)$ and to $(0,1]$ are $\alpha$-H\"older
continuous for some $0< \alpha \leq 1$, then
the derivative
$\tilde{f_{\ep}}^{\prime }$ is $\alpha$-H\"older continuous.
The $\alpha$-H\"older constant of $\tilde{f}_{\ep}$ depends on $\ep$
and may go to infinity as $\ep$ goes to zero.

\vskip5pt
\noindent{\sc Lemma 8.} {\em Suppose $\tilde{f}_{\ep}$ is a
continuous
mapping from $[-1,1]$ to the real line with a unique turning point
$0$.  Suppose
$\tilde{f}_{\ep}$ maps $1$
and $-1$ to $-1$ and $\tilde{f}_{\ep}(0) =1+\ep $.
If the restrictions of $\tilde{f}_{\ep}$ to $[-1,0]$ and to $[0,1]$ are
$C^{1}$ embeddings, then
$f_{\ep}=h_{\gamma, \ep}^{-1}\circ \tilde{f}_{\ep }\circ h_{\gamma ,
\ep }$ for any $\gamma >1$ is
a $C^{1}$ mapping from $[-1,1]$ to the real line and has the power law
$|x|^{\gamma }$ at the critical point.}

\vskip5pt
The proofs of Lemma 7 and Lemma 8 are the same as
those of Lemma 1 and 2.

\vskip5pt
\nd {\sc Lemma 9.} {\em Suppose $f_{\ep}$ has the power law
$|x|^{\gamma }$ with $\gamma >1$.  The
restrictions of $\tilde{f}_{\ep}$ to $[-1, 0]$ and $[0, 1]$ are
$C^{1+\alpha }$ embeddings for some $0< \alpha \leq 1$ if and
only if
the restrictions of $r_{\ep}$ to $[-1,0)$ and to $(0, 1]$
are $\alpha'$-H\"older continuous
for some $0< \alpha' \leq 1$.}

\vskip7pt
The proof of Lemma 9 is similar to Remark 1 and Remark 2.

\vskip10pt
\nd {\bf \S 4.1 Good families of mappings in ${\cal BH}\cup {\cal H}$}

\vskip5pt
Suppose $\{ f_{\ep}\}_{0\leq \ep \leq \ep_{0}}$ is a
family in ${\cal H}\cup {\cal BH}$ where $f_{\ep }(0)
=1+\ep$.

\vskip7pt
\nd {\sc Definition 7.} {\em
The family
$\{ f_{\ep}\}_{0\leq \ep \leq \ep_{0}}$ is a
\underline {good family} if it satisfies
the following conditions:

\begin{enumerate}
\item the mapping $F(x, \ep )=f_{\ep }(x)$ is $C^{1}$ in both variables
$x$ in $[-1,1]$ and $\ep$ in $[0, \ep_{0}] $,

\item each $ f_{\ep}$ has the same power law
$|x|^{\gamma }$ with $\gamma >1$ at
the critical points and
the functions
$R^{-}(x, \ep)=r_{\ep}(x)$ defined on $[-1,0]\times [0,\ep_{0}]$ and
$R^{+}(x, \ep)=r_{\ep}(x)$ defined on $[0,1]\times [0,\ep_{0}]$ are
continuous,

\item there are positive constants $K'$ and $\alpha' \leq 1$
such that $f_{\ep }$ is
$C^{1+\alpha' }$ and
the $\alpha'$-H\"older
constant of $f'_{\ep}$ is less than $K'$ for any
$0\leq \ep \leq \ep_{0}$,

\item there are positive constants $K''$ and $\alpha'' \leq
1$
such that the restrictions of $r_{\ep}$ to $[-1,0)$ and to $(0,1]$ are
$\alpha''$-H\"older continuous
and the
$\alpha''$-H\"older constants of these restrictions
are less than $K''$ for any $0\leq \ep \leq 1$,

\item
there are two \pocs\ $C_{0}$ and $\lambda <1$ such that
$\lambda_{n,\ep}\leq
C_{0}\lambda^{n}$ for all positive integers $n$ and $0\leq \ep \leq
\ep_{0}$.
\end{enumerate}

\nd Let $\alpha $ be the
minimum of $\alpha'$ and $\alpha''$.}

\vskip5pt
An example of a good family in ${\cal BH}\cup {\cal H}$ follows the
following proposition.

\vskip7pt
\nd {\sc Proposition 1.}
{\em Assume that
(a) $\{ F_{\ep}\}_{0\leq \ep \leq \ep_{0}}$ is a family of
$C^{3}$ embeddings on
$[-1,0]$ with nonpositive Schwarzian derivatives,
(b) $F_{\ep }$ fixes $-1$
and maps $0$ to $1+\ep$ for any $0\leq \ep \leq \ep_{0}$,
(c)
the derivative $F_{\ep }'(-1)$ of $F_{\ep}$ at $-1$ is greater than
$1/\gamma
$ and
(d) $G(x,\ep)=F_{\ep}(x)$ is $C^{2}$ on $[-1,1]\times
[0,\ep_{0}]$.
If $f_{\ep}(x)=F_{\ep}(-|x|^{\gamma })$,
then the family $\{ F_{\ep}(-|x|^{\gamma })
\}_{0\leq \ep \leq \ep_{0}}$ is a good family.}

\vskip5pt
To prove this proposition,  we only need to check the condition
(5) in Definition 7.
The condition (5) is a direct consequence of the following lemmas.

\vskip5pt
The first lemma is the $C^{3}$-Koebe distortion lemma. Suppose $I$ and $J$ 
are two intervals and $g$ is a
$C^{3}$ diffeomorphism from $I$ to $J$.
A measure of the nonlinearity of
$g$ is the function $n(g) =g^{\prime \prime }/g^{\prime }$.
If the absolute value of $n(g)$ on $I$
is bounded above
by a positive constant $C$,
then the distortion $|g'(x)|/|g'(y)|$ of $g$
at any pair
$x$ and $y$ in $I$ is bounded above by $exp(C|x-y|)$.
Suppose $d_{I}(x)$
is the distance from $x$ to the
boundary of $I$.

\vskip5pt
\noindent {\sc Lemma 10} (the $C^{3}$ Koebe distortion lemma). {\em
			  Suppose $g$ has nonnegative Scharzian derivative. Then
$n(g)(x)$
is bounded above by $2/d_{I}(x)$ for any $x$ in $I$.}

\vskip5pt
\noindent {\sl Proof.} See, for example, [J1] or [J5].

\vskip7pt
\nd {\sc Lemma 11.} {\em Suppose $\{ f_{\ep }\}_{0\leq \ep \leq
\ep_{0}}$ is the family in Proposition 1 and 
$\{ \eta_{n,
\ep}\}_{n=0}^{\infty }$ is the sequence determined by $f_{\ep}$ for
every $0\leq \ep \leq \ep_{0}$.
There
is a \poc\ $C$ which does not depend on parameter $\ep$ such that
for any $0\leq \ep \leq \ep_{0}$ and any pair $(J, I)$ with $J\subset
I$, $J\in \eta_{n+1, \ep}$ and $I\in \eta_{n, \ep}$, $|J|/|I| \geq C$.}

\vskip5pt
\nd {\sl Proof.} We suppress $\ep $ if there can be
no confusion. For any $0\leq \ep \leq \ep_{0}$, the first
partition
$\eta_{1}$ contains four intervals $I_{ 00}$, $I_{
01}$, $I_{11}$ and $I_{10}$.  There is a positive constant
$C_{1}$ which does not depend on $\ep$ such that the lengths of the left
interval $I_{00}$ and the
right interval $I_{10}$ are greater than $C_{1}$.
The $C^{3}$-Koebe distortion lemma says
$n(g_{
w})(x) \leq 2/d_{[-1,1]}(x)$
for any finite string $w$ of zeroes and ones. Moreover,
$n(g_{
w})(x) \leq 2/C_{1}$ if $x$ is in the union of two middle intervals
$I_{01}$ and $I_{11}$.
We also can find a constant $\tau >1$ which does not depend on $\ep$
such that $|f_{\ep }'(x)| \geq \tau $ for
all
$x$ in the union of the left interval $I_{00}$ and the right
interval
$I_{10}$.  Now the proof just follows the proof of Example 1 in [J2]. \Q

\vskip10pt
\nd {\bf \S 3.2 Asymptotic scaling function geometry of Cantor sets}

\vskip5pt
Let ${\cal A}$ stand the countable set of points in ${\cal C}^{*}$ whose
coordinates are
eventually all zeroes and let ${\cal B}$ stand the
complement of ${\cal A}$ in ${\cal C}^{*}$.

\vskip10pt
\noindent {\sc Theorem A.} {\em
Suppose $\{ f_{\ep}\}_{0\leq \ep \leq \ep_{0}}$ is a good family.
There is a family of H\"older
continuous functions $\{ s_{\varepsilon }
\}_{0< \ep \leq \ep_{0}}$
on the dual Cantor set ${\cal C}^{*}$
such that $s_{\varepsilon }$ is the
scaling function of $f_{\varepsilon }$ for any $0< \ep_{0}\leq \ep_{0}$,
and

\begin{enumerate}

\item for every $0< \ep_{1} \leq \ep_{0}$,
$s_{\ep }$ converges to
$s_{\ep_{1}}$ uniformly on ${\cal C}^{*}$ as $\ep $ tends to $\ep_{1}$,

\item for every $a^{*}\in {\cal C}^{*}$, the limit $s_{0}(a^{*})$ of $\{
s_{\ep}(a^{*})\}_{0< \ep \leq \ep_{0}}$ as $\ep$ decreases to zero
exists, the limiting function $s_{0}(a^{*})$
is the scaling function of $f_{0}$ and satisfies:

\begin{description}

\item[{\em 2.1.}] $s_{0}$ has jump discontinuities at all points in
${\cal A}$,

\item[{\em 2.2.}] $s_{0}$ is continuous at all points in ${\cal B}$ and
the restriction of $s_{0}$ to ${\cal B}$ is a H\"older continuous
function.
\end{description}
\end{enumerate} }

We will prove this theorem through several lemmas. The first lemma is
similar to the $C^{1+\alpha}$-Denjoy-Koebe distortion lemma in [J1].
We call it
the \underline {uniform $C^{1+\alpha}$-Denjoy-Koebe distortion
lemma}. 

\vskip5pt
Suppose $\{ f_{\ep} \}_{0\leq \ep \leq \ep_{0}}$ is a family of
mappings in ${\cal H}\cup {\cal BH}$ and satisfies the
\tcs\ (1)--(4) in Definition 7.  We suppress $\ep$ when
there can be no confusion.
For each $0\leq \ep \leq \ep_{0}$,  $\eta_{1}$
contains four
intervals $I_{00}$, $I_{01}$, $I_{11}$ and $I_{
10}$. Suppose $x$ and $y$ are in one of these four intervals
and $J_{0}$ is the interval bounded by $x$ and $y$.
Let
$\theta(x,y) =\{ J_{0}, J_{1},\cdots \}$ be a
sequence of backward images of $J_{0}$ under
$f_{\ep}$, that
means, the restriction of $f_{\ep }$ to $J_{n}$ embeds $J_{
n}$ onto $J_{n-1}$ for any positive integer $n$.
Let $g_{n}$ be the inverse
of the restriction of the $n^{th}$ iterate of $f_{\ep }$ to
$J_{n}$.
Let $d_{xy}$ be the distance from $\{ x,y\}$ to $\{ -1, 1\}$.
Define the distortion of the $n^{th}$ iterate of
$f_{\ep }$ at $x$ and $y$ along
$\theta(x, y)$ to be the ratio $|g_{n}^{\prime
}(x)|/|g_{n}^{\prime }(y)|$.

\vskip7pt
\nd {\sc Lemma 12.} (the uniform $C^{1+\alpha}$-Denjoy-Koebe distortion lemma.) {\em There are positive constants $A$, $B$ and $C$ such
that for any $0\leq \ep \leq \ep_{0}$,	any  $x$ and $y$ in one
of the intervals in $\eta_{\ep ,1}$
and any sequence of backward images
$\theta_{\ep}(x, y)=\{ J_{\ep , 0}, J_{\ep , 1}, \cdots \}$ of $J_{\ep ,
0}$ under $f_{\ep}$, the distortion of the $n^{th}$ iterate of
$f_{\ep}$ at $x$ and $y$ along
$\theta_{\ep }(x, y) $ satisfies
\[ \frac{|g_{\ep , n}^{\prime }(x)|}{|g_{\ep ,n}^{\prime }(y)|}
\leq
\exp ((A +B\sum_{i=1}^{n}|J_{\ep , i}|+
\frac{C|J_{\ep , 0}|}{d_{xy}})
\sum_{i=1}^{n}|J_{\ep , i}|^{\alpha})\]
for every \poi\ $n$.}

\vskip7pt
\nd {\sl Proof.}
For every $0\leq \ep \leq \ep_{0}$, $\eta_{2}$
contains
eight intervals. The two of them which close to $0$ are $I_{010}$ and
$I_{110}$.  Suppose $
I_{010}= [a, b]$ and
$I_{110}=[c, d]$.  We call $[a, 0]$ and $[0,
d]$ the middle intervals,  $[-1, a]$ the left interval
and $[d, 1]$ the right interval
(see Figure 7).

\vskip5pt
\begin{picture}(300,70)(0,0)
\put(50,50){\line(1,0){25}}
\put(80,50){\line(1,0){25}}
\put(110,50){\line(1,0){25}}
\put(140,50){\line(1,0){25}}
\put(170,50){\line(1,0){25}}
\put(200,50){\line(1,0){25}}
\put(230,50){\line(1,0){25}}
\put(260,50){\line(1,0){25}}

\put(50,55){$I_{\ep ,000}$}
\put(80,55){$I_{\ep ,001}$}
\put(110,55){$I_{\ep ,011}$}
\put(140,55){$I_{\ep ,010}$}
\put(170,55){$I_{\ep ,110}$}
\put(200,55){$I_{\ep ,111}$}
\put(230,55){$I_{\ep ,101}$}
\put(260,55){$I_{\ep ,100}$}
\put(167,50){$\cdot$}
\put(167,42){$0$}

\put(80,20){left}
\put(150,20){middle}
\put(230,20){right}

\put(130,0){Figure 7}
\end{picture}

\vskip10pt
By the condition (1) in Definition 7, there is a positive
constant $C_{1}$ which does not depend on $\ep$ such that
the lengths of the left interval and the right interval
are greater than $C_{1}$.

\vskip10pt
By the
conditions (1) and (3) in Definition 7, there are positives constants
$c_{1}$ and $K_{1}$ which do not depend on $\ep$
such that the minimum value of $|f_{\ep}'|$ on the union of the
left and right intervals is greater than $c_{1}$ and the
$\alpha'$-Holder constants of the restrictions of $f_{\ep}'$ to the left
interval and to the right interval are less than $K_{1}$.

\vskip10pt
Suppose $y$ is not one of $0$, $1$ and $-1$ and
$x$ is the preimage of $y$ under $h_{\gamma , \ep}$. By the chain rule,

\[ \tilde{f}_{\ep}^{\prime }(y) = \frac{ f_{\ep}^{\prime }(x)
((1+\ep)^{2} -x^{2})^{\frac{\gamma-1}{\gamma}}}{((1+\ep)^{2} -
(f_{\ep}(x))^{2})^{\frac{\gamma-1}{\gamma}}}.\]
This equation and the conditions $(1)$ and $(4)$ in
Definition 7 imply the restrictions of
$\tilde{f}_{\ep}$ to $[-1, 0]$ and to $[0, 1]$ are
$C^{1+\alpha '' }$ embeddings. Moreover,
there
are constants $c_{2}$ and $K_{2}$ which do not depend on $\ep$ such
that the minimum value
of $|\tilde{f}^{\prime }_{\ep}|$ on the image of every one of the middle
intervals under $h_{\gamma , \ep}$
is greater then $c_{2}$ and the $\alpha''$-H\"older constant of
the restriction of $\tilde{f}_{\ep}' $ to the image of every one of the
middle intervals is less than $K_{2}$.

\vskip10pt
The restriction of $h_{\gamma , \ep}^{\prime }$ to the union of the
middle intervals is Lipschitz continuous.  There are positive constants
$c_{3}$ and $K_{3}$ which do not depend on $\ep$ such that the minimum
value of restriction of
$h_{\gamma , \ep}^{\prime }$ to the union of middle intervals is greater
then $c_{3}$ and the Lipschitz constant of such restriction is less than
$K_{3}$.

\vskip5pt
Let $x_{i}$ and $y_{i}$ be
the images of $x$ and $y$ under $g_{\ep ,i}$.
Notice that this implies that $x_{i}$ and
$y_{i}$ lie in the same interval of $\eta_{i+1}$ for
$i\geq 0$.   For every integer $n>0$,  $g_{\ep ,n}'(x)/g_{\ep ,n}'(y)$
equals $(f_{\ep}^{\circ n})'(y_{n})/(f_{\ep}^{\circ n})'(x_{n})$. By
the chain
rule, the ratio $(f_{\ep}^{\circ n})'(y_{n})/(f_{\ep}^{\circ
n})'(x_{n})$
equals the product of ratios $f_{\ep}'(y_{n-i})/f_{\ep}'(x_{n-i})$ where
$i$ runs from $0$ to $n-1$.   This product can be factored into two
products,
\[ \prod_{x_{i}, y_{i} \in LR} f_{\ep}'(y_{i})/f_{\ep}'(x_{i}) \hskip7pt
and \hskip7pt \prod_{x_{i}, y_{i} \in M}
f_{\ep}'(y_{i})/f_{\ep}'(x_{i}).\]
Here $LR$ stands for the union of the left and right
intervals and
$M$ stands for the union of the two middle intervals.
We factor the product
$\prod_{x_{i}, y_{i} \in M} f_{\ep}'(y_{i})/f_{\ep}'(x_{i})$
into three factors,
\[ \prod_{x_{i}, y_{i} \in M} \frac{f_{\ep}'(y_{i})}{f_{\ep}'(x_{i})}=
\prod_{x_{i}, y_{i}\in M}\frac{h_{\gamma , \ep
}^{\prime }(y_{i})}
{h^{\prime }_{\gamma , \ep}(x_{i})}\cdot
\prod_{x_{i}, y_{i}\in M}\frac{\tilde{f}_{\ep}'(h_{\gamma ,\ep}(y_{i}))}
{\tilde{f}_{\ep}'(h_{\gamma , \ep}(x_{i}))}\cdot
\prod_{x_{i},y_{i}\in M}\frac{h^{\prime }_{\gamma ,\ep}(f(x_{i}))}
{h^{\prime}_{\gamma , \ep}(f(y_{i}))}.\]
The third factor of them
can be factored again into two
products,
\[ \prod_{x_{i}, y_{i}\in M}
\frac{(1+f_{\ep}(y_{i}))^{\frac{\gamma-1}{\gamma}
}}{(1+f_{\ep}(x_{i}))^{\frac{\gamma-1}{\gamma
}}}\hskip8pt and \hskip8pt \prod_{x_{i}, y_{i}\in M}
\frac{(1-f_{\ep}(y_{i}))^{\frac{\gamma-1}{\gamma
}}}{(1-f_{\ep}(x_{i}))^{\frac{\gamma-1}{\gamma }}}.\]

\vskip10pt
Now just following the arguments in the proof of the
$C^{1+\alpha}$-Denjoy-Koebe distortion lemma in [J1],
we can estimate every factors.	We then put all estimations
together to get $A=
K_{1}/c_{1}+(K_{3}^{\alpha }K_{2})/c_{2}
+K_{3}/c_{3} +(\gamma-1)/\gamma$,
$B=
(\gamma-1)/(\gamma C_{1})$ and $C=
(\gamma-1)/\gamma$. The constants $A$, $B$ and $C$ do not depend on the
parameter $\ep $. \Q

\vskip10pt
\nd {\sc Corollary 1.} {\em If the family $\{ f_{\ep}\}_{0\leq \ep
\leq
\ep_{0}}$ is a good family,
then there are positive constants $D$ and
$E$ such that for any $0\leq \ep \leq \ep_{0}$,
any
$x$ and $y$ in one of the intervals in $\eta_{1, \ep}$
and any sequence $\theta_{\ep}(x, y)=\{ J_{\ep , 0}, J_{\ep ,
1},\cdots \} $ of backward images of $J_{\ep , 0}$ under $f_{\ep}$, the
distortion of the $n^{th}$ iterate of $f_{\ep}$ along
$\theta_{\ep} (x,y)$ satisfies
\[ \frac{|g_{\ep , n}^{\prime }(x)|}{|g_{\ep ,n}^{\prime }(y)|}\leq
\exp (D
+ \frac{E}{d_{x
y}})|J_{\ep , 0}|^{\alpha }\]
for every \poi\ $n$.}

\vskip10pt
\nd {\sl Proof.}
The condition (5) in Definition 7 implies $\sum_{i=0}^{n} |J_{\ep ,
i}|$ is less that $C_{2}= 2 C_{0}/(1-\lambda)$ and
$\sum_{i=0}^{n} |J_{\ep , i}|^{\alpha }$ is less that
$P= 2^{\alpha }|J_{\ep , 0}|^{\alpha }C_{0}/(1-\lambda^{\alpha
})=
|J_{\ep , 0}|^{\alpha }C_{3}$.
We use $D$ to denote
$(A+ BC_{2})C_{3}$ and $E$ to denote $CC_{3}$, where $A$, $B$ and $C$
are the constants in the previous lemma.  Now it is easy to show
this corollary from the previous
lemma. \Q

\vskip10pt
Before we prove more lemmas, we study asymptotic behavior of the
maximal invariant set $\Lambda_{\ep}$ of $f_{\ep}$ in ${\cal H}$ when $f_{\ep}$
approaches ${\cal BH}$.

\vskip10pt
\nd {\bf \S 4.3 Determination of the geometry of Cantor set by the
leading gap.}

\vskip5pt
Suppose $f_{\ep}$
is a mapping in ${\cal H}$. We suppress $\ep$ when there can be
no confusion. Let $\Lambda$
be the maximal invariant set of $f_{\ep}$ and $\{
\eta_{n
}\}_{n=0}^{\infty}$ be the sequence determined by $f_{\ep}$.
For any positive integer $n$ and any
$I_{w}$
in $\eta_{n}$,	let $I_{w0}$ and $I_{w1}$ be the
two
intervals in $\eta_{n+1}$ which are contained in $I_{w}$.  We
call the complement of $I_{w0}$ and $I_{
w1}$ in $I_{w}$ the \underline {gap} on $I_{w}$ and denote it by $G_{
 w}$.	Let $G$ be the complement of $I_{0}$ and
$I_{1}$ in $[-1, 1]$. We call $G$ the \underline {leading gap}.

\vskip5pt
\nd {\sc Definition 8.} {\em We call the set of ratios, $\{ |G_{\ep ,
w}|/|I_{\ep ,w}|\}$, for all finite strings $w$
of zeroes and ones the \underline {gap geometry} of $\Lambda_{\ep}$ or
$f_{\ep}$.}

\vskip5pt
Suppose $\{ f_{\ep}\}_{0\leq \ep \leq \ep_{0}}$ is a family
in ${\cal H} \cup {\cal BH}$
and $\{ \Lambda_{\ep} \}_{0\leq \ep \leq
\ep_{0}}$ is the family of the corresponding maximum invariant sets.

\vskip5pt
\nd {\sc Definition 9.} {\em Suppose $\beta $ is a function defined on
$[0, \ep_{0}]$.  We say $\beta$ determines asymptotically the gap
geometry of $\{ \Lambda_{\ep}\}_{0<
\ep \leq 1}$
if there is a positive constant $C$ such that
for all
$0\leq \ep \leq \ep_{0}$, all finite strings $w$ of zeroes and
ones and $i=0$ or $1$,

\begin{enumerate}
\item $C^{-1}\beta (\ep )\leq
|G_{\ep , w}|/|I_{\ep ,w}|	 \leq C\beta (\ep )$,

\item $|I_{\ep , wi}|/|I_{\ep ,w}|\geq C^{-1}$.
\end{enumerate}
\nd The constant $C$ is called a determining constant.}

\vskip5pt
\nd {\sc Theorem B.} {\em
Suppose $\{ f_{\ep}\}_{0\leq \ep \leq \ep_{0}}$ is a good
family.  Then
the family $\{
\Lambda_{\ep}\}_{0< \ep \leq \ep_{0}}$ is a family of Cantor sets.
Moreover, the function
$\ep^{\frac{1}{\gamma}}$ on $[0,\ep_{0}]$
determines asymptotically the gap geometry of
$\{ \Lambda_{\ep}\}_{0\leq
\ep \leq \ep_{0}}$.}

\vskip5pt
\nd {\sl Proof.}
For each $0\leq \ep \leq \ep_{0}$,  $\eta_{1}$ contains four
interval $I_{00}$, $I_{01}$, $I_{11}$ and $I_{
10}$.  We call $I_{00}$ the left interval,
$I_{10}$ the right interval and
$I_{01}$ and
$I_{11}$ the middle intervals. We also call $\eta_{1}$ the
first level and $\eta_{2}$ the second level.

\vskip5pt
By the \tcs\ (2) in Definition 7,
there is a \poc\ $C_{1}$ which does not depend on $\ep$ such that
$C_{1}^{-1}\ep^{\frac{1}{\gamma}} \leq |G_{0}|\leq
C_{1}\ep^{\frac{1}{\gamma}}$.

\vskip5pt
By \tc\ (1) in Definition 7, there is a \poc\ $C_{2}$ which does not depend on $\ep$
such that for any triple $(G, J, I)$ where $I$ is an interval in the
first level or $\eta_{0}$, $G$ is the gap on $I$ and $J\subset I$ is an
interval in the second interval or the first level, $|J|/|I|\geq C_{2}$
and $C_{2}^{-1}\ep^{\frac{1}{\gamma}}
\leq |G|/
|I|\leq C_{2}\ep^{\frac{1}{\gamma}}$.

\vskip5pt
For any integer
$n>1$ and any triple
$(G, J,
I)$ where $I$ is in $\eta_{n}$,
$G$ is the gap on $I$ and $J\subset I$ is in
$\eta_{n+1}$,  let $G_{i}$, $J_{i}$ and $I_{i}$ be the images of
$G$, $J$ and $I$ under the $i^{th}$ iterate of $f_{\ep}$ for
$0\leq i\leq n-1$.  Then $I_{n-1}$ is in the
first level, $G_{n-1}$ is the gap on $I_{n-1}$ and $J_{n-1}\subset
I_{n-1}$ is in the second level.

\vskip5pt
We divide the possible itineraries of the sequence of
triples $\{ (G_{i}$, $ J_{i}$, $I_{i})
\}_{i=0}^{n-1}$ into two cases.  The first case is that no one of
$I_{i}$
is in the union of the middle intervals. The second case is that some of
$I_{i}$ is in one of the middle intervals.

\vskip5pt
In the first case, $G_{i}$, $J_{i}$ and $I_{i}$ are in the union of the
left interval and the right interval for all $0\leq i<n$.
By the \tcs\ (1) and (3) in Definition 7, there is a positive constant
$c$ which does not depend on $\ep$ such that
the minimum value of the restriction of $f_{\ep}'$ to the
union of the left and right intervals is greater than
$c$.
Suppose $\lambda $, $C_{0}$, $\alpha'$ and $K'$ are
the constants in the \tcs\ $(3)$ and $(5)$ in Definition 7. By the
naive distortion lemma (see [J1]), there is a constant
$C_{3}$ which equals
$C_{2}\exp(K'C_{0}/(c
(1-\lambda^{\alpha'})))$ such that
$|J|/|I|\geq C_{3}^{-1}$ and
$C_{3}^{-1}\ep^{\frac{1}{\gamma}}
\leq |G|/
|I|\leq C_{3}\ep^{\frac{1}{\gamma}}$ because
$G$, $J$ and $I$ are the images of
$G_{n-1}$, $J_{n-1}$ and $I_{n-1}$ under the $(n-1)^{th}$ iterate of
$f_{\ep}$.

\vskip5pt
In the second case, let $m$ be the largest positive integer such that
$I_{m}$ is in one of the middle intervals.  We can divide this case into
two subcases according to $m$.	One is that $m$ is $n-1$. The other is
that $m$ is less than $n-1$.

\vskip5pt
If $m=n-1$, then $I_{n-1}$ is one of the middle intervals.  By the \tc\
(1) in Definition 7, there is a positive constant $C_{4}$ which does
not depend on $\ep$ such that the lengths
of the left interval and the right interval are greater then $C_{4}$.
By Corollary 1,
there is a constant $C_{5}$ which equals $C_{2}exp (D+E/C_{4})$ and does
not depend on $\ep$ such that $|J|/|I|\geq C_{5}^{-1}$ and
$C_{5}^{-1}\ep^{\frac{1}{\gamma}}
\leq |G|/
|I|\leq C_{5}\ep^{\frac{1}{\gamma}}$ because the restriction of the
$(n-1)^{th}$ iterate of $f$ to $I$ embeds $I$ to $I_{n-1}$
and the distance from $I_{n-1}$ to $\{ -1, 1\} $ is greater than
$C_{4}$.

\vskip5pt
If $m< n-1$, then $I_{i}$ is in the union of the left and right
intervals for $m< i\leq n-1$.  Because $I_{n-1}$ is one of the left and
the right intervals,
$I_{m+1}$ has $1$ as a boundary and
$I_{m}$ is the one closing
$0$ in $\eta_{n-m, \ep }$.
For the sequence $I_{m+1}$, $\cdots $, $I_{n-1}$, no one of them
is in the union of the middle intervals. By the same arguments as
those in the first case imply that
$|J_{m+1}|/|I_{m+1}|\geq C_{3}^{-1}$ and
$C_{3}^{-1}\ep^{\frac{1}{\gamma}}
\leq |G_{m+1}|/
|I_{m+1}|\leq C_{3}\ep^{\frac{1}{\gamma}}$.  For the sequence $I_{0}$,
$\cdots
$, $I_{m}$, the last one $I_{m}$ is in one of the middle intervals.
Similar arguments to those in the subcase $m=n-1$ imply that
$|J|/|I|\geq C_{6}^{-1}|J_{m}|/|I_{m}|$ and
$C_{6}^{-1}|G_{m}|/|I_{m}|
\leq |G|/
|I|\leq C_{6}|G_{m}|/|I_{m}|$, where
$C_{6} = exp (D+E/C_{4})$.
Because the restriction of $f_{\ep}$ to
$I_{m}$ is comparable with
the power law mapping $|x|^{\gamma }$ uniformly on $\ep $ by the \tc\
(2) in Definition 7, we may assume
$f_{\ep}|I_{m} =1+\ep -|x|^{\gamma }$.	There is a positive constant $C_{7}$ which does not
depend on $\ep$ such that
$|J_{m}|/|I_{m}|\geq C_{7}^{-1}|J_{m+1}|/|I_{m+1}|$ and
$C_{7}^{-1}|G_{m+1}|/|I_{m+1}|
\leq |G_{m}|/
|I_{m}|\leq C_{7}|G_{m+1}|/|I_{m+1}|$.
We get that
$|J|/|I|\geq (C_{6}C_{7}C_{3})^{-1}$ and
$(C_{6}C_{7}C_{3})^{-1}\ep^{\frac{1}{\gamma}}
\leq |G|/
|I|\leq (C_{6}C_{7}C_{3})\ep^{\frac{1}{\gamma}}$.

\vskip5pt
Let $C$ be $C_{7}C_{6}C_{3}$. It is a determining constant of the gap
geometry of
$\{\Lambda_{\ep }\}_{0\leq \ep \leq \ep_{0}}$. \Q

\vskip5pt
Let
$HD(\ep)$ be the Hausdorff dimension
of $\Lambda_{\ep}$.  An immediate consequence of Theorem B is
the following corollary.

\vskip5pt
\nd {\sc Corollary 2.} {\em
There is a \poc\ $C$ which does not depend on $\ep$ such that
\[ 0<  HD(\ep)\leq 1-C\ep^{\frac{1}{\gamma}}\]
for all
$0\leq \ep \leq \ep_{0}$.}

\vskip5pt
\nd {\sc Corollary 3.}
{\em There is a \poc\ $C$ which does not depend on $\ep$ such
that
\[ C^{-1}\ep^{\frac{1}{\gamma}} \leq
s_{\ep}((a^{*}0.))+s_{\ep}((a^{*}1.)) \leq C\ep^{\frac{1}{\gamma}}\]
for all $a^{*}$ in ${\cal C}^{*}$ and $0\leq \ep \leq
\ep_{0}$.}

\vskip5pt
\nd {\sl Proof.} For any $0< \ep \leq \ep_{0}$ and any $a^{*} \in
{\cal C}^{*}$, $s(w1)+s(w0)= |G_{w}|/|I_{w}|$ where $a^{*}= (\cdots
w.)$ and $G_{w}$ is the gap on $I_{w}$.	Now this corollary is a
consequence of Theorem A and Lemma 5. \Q

\vskip10pt
\nd {\bf \S 4.4 The proof of Theorem A}

\vskip5pt
Let us go on to prove Theorem A.
Suppose $\{ f_{\ep} \}_{0\leq \ep \leq \ep_{0}}$ is a good family.
For any
$a^{*}$
in ${\cal C}^{*}$, we always use $w_{n}i$ to denote the
first $(n+1)$ coordinates
of $a^{*}$, that is, $a^{*}=(\cdots w_{n}i.)$, and use $s(\ep ,
w_{n}i)$ to denote the ratio of lengths, $|I_{\ep , w_{n}i}|/
|I_{\ep , w_{n}}|$.  Suppose $s_{\ep}$ is the
scaling function of $f_{\ep}$ for $0< \ep \leq \ep_{0}$. We suppress
$\ep$ when there can be no confusion.

\vskip5pt
\nd {\sc Lemma A1.} {\em For each $0< \ep_{1} \leq \ep_{0}$,
$s_{\ep}$ converges to $s_{\ep_{1}}$ uniformly on ${\cal C}^{*}$ as $\ep $
tends to $\ep_{1}$.}

\vskip5pt
\nd {\sl Proof.} Suppose $C_{\ep}$ is the constant
obtained in Lemma 5.	By the \tc\ (1) in Definition 7,
$C_{\ep}$ is continuous on $0< \ep \leq \ep_{0}$.  We can find a
positive number
$\delta $ such that $C_{\ep }\leq 2C_{\ep_{1}}$ for all $\ep$ in
$(\ep_{1}-\delta , \ep_{1}+\delta)$.

\vskip5pt
From the proof of Lemma 5, we have that for any $a^{*}$ in
${\cal C}^{*}$, $|s_{\ep}(a^{*})-s(\ep, w_{n}i)|
\leq C_{\ep} |I_{\ep , w_{n}}|^{\alpha' }$ and
$|s_{\ep_{1}}(a^{*})-s(\ep_{1}, w_{n}i)|
\leq C_{\ep_{1}} |I_{\ep_{1} , w_{n}}|^{\alpha' }$.
Because we may write
$|s_{\ep}(a^{*})-s_{\ep_{1}}(a^{*})|$ in
$|s_{\ep}(a^{*})-s(\ep , w_{n}i) +s(\ep ,w_{n}i) -s(\ep_{1}, w_{n}i)
+s(\ep_{1}, w_{n}i) - s_{\ep_{1}}(a^{*})|$, this implies that
\[ |s_{\ep}(a^{*})-s_{\ep_{1}}(a^{*})|\leq
2C_{\ep_{1}}|I_{\ep , w_{n}}|^{\alpha' }+C_{\ep_{1}} |I_{\ep_{1} ,
w_{n}}|^{\alpha' }+
|s(\ep , w_{n}i)-s(\ep_{1}, w_{n}i)|\]
for all $n >0$ and $\ep$ in $(\ep_{1}-\delta , \ep_{1}-\delta )$.  Now
the last inequality and the \tcs\ $(1)$ and $(5)$ in Definition 7
imply this lemma. \Q

\vskip5pt
\nd {\sc Lemma A2.} {\em For every $a^{*}$ in ${\cal
C}^{*}$, the limit
of $\{ s_{\ep}(a^{*}) \}_{0<\ep \leq \ep_{0}}$ exists as $\ep$ decreases
to zero.}

\vskip5pt
\nd {\sl Proof.} For every $0\leq \ep \leq \ep_{0}$,
$\eta_{1, \ep }$
contains four intervals $I_{\ep , 00}$, $I_{\ep ,01}$, $I_{\ep ,11}$
and $I_{\ep ,10}$.  We call
$I_{\ep ,01}$ and  $I_{\ep ,11}$ the middle intervals,
$I_{\ep ,00}$ the left interval and
$I_{\ep ,10}$ the right interval.

\vskip5pt
For $a^{*}$ in $C^{*}$, we may arrange it into two cases according to
its coordinates. The first case is that
the coordinates of $a^{*}$ are
eventually all zeroes. The second case is that there are infinite many
ones in the coordinates of $a^{*}$.

\vskip5pt
In the first case, we can find a \poi\ $N$ such that
$I_{\ep , w_{n}}$ is in the left interval
for every $n\geq N$. By the \tc\ (1) in Definition 7,  there is a
positive
constant $c$ which does not depend on $\ep$ such that the minimum value
of $|f_{\ep}'|$ on the left interval is greater than $c$.
By the naive distortion lemma and and the \tcs\ $(3)$ and $(5)$ in
Definition 7 and similar arguments to the proof of Lemma A1,  there
is a constant $C_{1}$, which equals
$exp( K'C_{0}/(c(1-\lambda^{\alpha' }))$, such that
\[ |s_{\ep}(a^{*})-s_{\ep'}(a^{*})|\leq
C_{1}(|I_{\ep , w_{n}}|^{\alpha' }+
|I_{\ep' , w_{n}}|^{\alpha' })+
|s(\ep , w_{n}i)-s(\ep' , w_{n}i)|\]
for all $\ep $ and $\ep'$ in $(0,\ep_{0}]$
and $n\geq N$. Now we can show that the limit of
$\{ s_{\ep}(a^{*}) \}_{0<\ep \leq \ep_{0}}$ as $\ep$ decreases to
zero exists.

\vskip5pt
In the second case, by the \tc\ $(1)$ in Definition 7,
there is a positive constant $C_{2}$ which does not depend on $\ep$ such
that
the lengths of the left interval and the right interval are greater than
$C_{2}$ for any $0\leq \ep \leq \ep_{0}$.  By Corollary 1, there is a
constant $C_{3}= D+E/C_{2}$ which does not depend on $\ep$ such that if
$I_{\ep , w_{n}}$ is in one of the middle
intervals, then $|s(\ep , w_{m}i)- s(\ep , w_{n}i)| \leq C_{3}|I_{\ep ,
w_{n}}|^{\alpha }$ for all $\ep$ in $(0, \ep_{0}]$ and $m>n>0$ because
$s(\ep , w_{m}i) = (|(f^{\circ (m-n)})'(x)|/
|(f^{\circ (m-n)})'(y)|)s(\ep , w_{n}i)$ for some $x$ and $y$ in $I_{\ep
, w_{m}}$.  Moreover, let $m$ increases
to infinity, then $|s_{\ep}(a^{*})- s(\ep , w_{n}i)| \leq C_{3}|I_{\ep ,
w_{n}}|^{\alpha }$
if $I_{\ep , w_{n}}$
is in one of the middle intervals.

\vskip5pt
If the $n^{th}$ coordinate of $a^{*}$ is one, then for any $\ep$ and
$\ep'$ in $(0, \ep_{0}]$,  $I_{\ep , w_{n}}$ is in one of the middle
intervals for $\ep$
and $I_{\ep' , w_{n}}$ is in one of the middle intervals for $\ep'$.
Then
\[ |s_{\ep}(a^{*})-s_{\ep'}(a^{*})|\leq
|s(\ep , w_{n}i)-s(\ep' , w_{n}i)|+
C_{3}(|I_{\ep , w_{n}}|^{\alpha}+
|I_{\ep' , w_{n}}|^{\alpha}).\]
Because there are infinite many ones in the coordinates of $a^{*}$,
the last inequality implies that the limit of $\{ s_{\ep}(a^{*})
\}_{0<\ep \leq \ep_{0}}$ as $\ep$ decreases to zero exists. \Q

\vskip5pt
Let $s_{0}(a^{*})$ be the limit of
$\{ s_{\ep}(a^{*}) \}_{0<\ep \leq \ep_{0}}$ as $\ep$ decreases to zero.
Then $s_{0}$ defines a function on
${\cal C}^{*}$.

\vskip5pt
\nd {\sc Lemma A3.} The limiting function $s_{0}$ is the scaling
function of $f_{0}$.

\vskip5pt
\nd {\sl Proof.} The proof is similar to the proof of Lemma A1.
Let us outline the proof. There are four intervals in the first
partition $\eta_{1}$ determined by $f_{0}$.  We call the one adjacent to
$-1$ the left interval,  the two adjacent to $0$
the middle intervals and the one adjacent to $1$ the
right interval. For any $a^{*}$ in ${\cal C}^{*}$, its
coordinates either are eventually all zeroes or contains infinite many
ones.

\vskip5pt
If the coordinates of $a^{*}$ are eventually all zeroes, then all $I_{
w_{n}}$ are eventually in
the left interval.  There is a positive constant $C_{1}$
such that for any $m>n>0$,
$|s(
w_{m}i)-s(w_{n}i)|\leq C_{1} |I_{w_{n}}|^{\alpha' }$ by the
naive distortion lemma.

\vskip5pt
If there are infinite many ones in the coordinates of $a^{*}$.
Suppose the $n^{th}$
coordinate of $a^{*}$ is one. Then
$I_{w_{n}}$
is in one of the middle interval.
There is positive constant $C_{2}$ such that for any $m>n>0$,
$|s(
w_{m}i)-s(w_{n}i)|\leq C_{2} |I_{w_{n}}|^{\alpha }$ by the
Corollary 1.

\vskip5pt
In both of cases, the limit of $\{ s(
w_{n}i)\}_{n=0}^{\infty}$ as the length of $w_{n}i$ increases to
infinity exists. The scaling function of $f_{0}$ exists.
Allowed $\ep=0$ in the proof of Lemma A2,
we can show that this scaling function is $s_{0}$. \Q

\vskip5pt
\nd {\sc Lemma A4.} {\em The scaling function $s_{0}$ has jump
discontinuities
at all points in ${\cal A}$.}

\vskip5pt
\nd {\sl Proof.} Suppose $a^{*}$ is in ${\cal A}$ and $a^{*} =
(0_{\infty}wi.)$ where $0_{\infty}$ is the one-sided infinite string
of zeroes extending to the left, $w$ is a finite string of zeroes and
ones and $i$ is either zero or one.  Let $0_{n}$ be the finite string of
zeroes of length $n$. The
interval $I_{0_{n}w}$ is eventually in the left interval $I_{00}$.
We use $b_{n}$ to
denote the length of $I_{0_{n}w}$ and $a_{n}$ or $a_{n}'$ to denote
the length of
$I_{0_{n}wi}$.	Let $c_{n}$ be the distance from $I_{
0_{n}w}$ to $-1$ (see Figure 8).

\vskip5pt
\begin{picture}(340,200)(0,0)
\put(10,175){$\cdot$}
\put(0,185){$-1$}
\put(30,180){\line(1,0){40}}
\put(40,187){$I_{0,0_{n}w}$}
\put(30,180){\line(0,-1){5}}
\put(50,180){\line(0,-1){5}}
\put(70,180){\line(0,-1){5}}
\put(33,160){$a_{n}$}
\put(53,160){$a^{'}_{n}$}
\put(13,160){$c_{n}$}
\put(45,145){$b_{n}$}

\put(340,125){$1$}
\put(340,117){$\cdot$}
\put(290,120){\line(1,0){40}}
\put(300,127){$I_{0,10_{n}w}$}
\put(290,120){\line(0,-1){5}}
\put(310,120){\line(0,-1){5}}
\put(330,120){\line(0,-1){5}}

\put(165,75){$0$}
\put(165,67){$\cdot$}
\put(115,70){\line(1,0){40}}
\put(120,77){$I_{0,010_{n}w}$}
\put(115,70){\line(0,-1){5}}
\put(135,70){\line(0,-1){5}}
\put(155,70){\line(0,-1){5}}

\put(175,70){\line(1,0){40}}
\put(180,77){$I_{0,110_{n}w}$}
\put(175,70){\line(0,-1){5}}
\put(195,70){\line(0,-1){5}}
\put(215,70){\line(0,-1){5}}

\put(285,125){\vector(-4, 1){200}}
\put(175,90){\vector(4, 1){90}}

\put(90,20){\vector(4, 3){45}}
\put(140,20){\vector(4, 3){45}}

\put(200,150){quasilinear}
\put(240,90){power law}
\put(180,40){quasilinear}
\put(140,0){Figure 8}

\end{picture}

\vskip5pt
Let $j$ be either zero or one. 	Because $f_{0}$ has the power law
$|x|^{\gamma}$ with $\gamma >1$ at the critical point and $I_{j10_{n}w}$
close to the critical point, the limit of $\{ s(
j10_{n}wi)\}_{n =0}^{\infty}$ equals the
limit of
\[ s_{n, 1}=
\frac{(a_{n}+c_{n})^{\frac{1}{\gamma}}-c_{n}^{\frac{1}{\gamma}}}
{(b_{n}+c_{n})^{\frac{1}{\gamma}}-c_{n}^{\frac{1}{\gamma}}}\]
or the limit of
\[
s_{n,2}=
\frac{(b_{n}+c_{n})^{\frac{1}{\gamma}}-(a_{n}+c_{n})^{\frac{1}{\gamma}}}
{(b_{n}+c_{n})^{\frac{1}{\gamma}}-c_{n}^{\frac{1}{\gamma}}}\]
as $n$ increases to infinity if the limits of $s_{n, 1}$ and
$s_{n,2}$ as $n$ increases infinity exist.

\vskip5pt
Because the minimum value of the restriction of $f_{0}$ to the left
interval $I_{00}$ is positive, by using the naive distortion
lemma, we can show the limit of $\{b_{n}/c_{n}\}_{n=0}^{\infty}$ and
the limit $\{ a_{n}/c_{n}\}_{n=0}^{\infty}$ as $n$ increases to
infinity exist.  We use $\tau_{1}(a^{*})$ and $\tau_{2}(a^{*})$ to
denote these limit, respectively. Now we conclude that
the limit of $\{ s(
j10_{n}wi)\}_{n =0}^{\infty}$ as $n$ increases to infinity exists and
\[\lim_{n\mapsto +\infty } s(j10_{n}wi)=
\frac{(1+\tau_{2}(a^{*}))^{\frac{1}{\gamma}}-1}
{(1+\tau_{1}(a^{*}))^{\frac{1}{\gamma}}-1} \hskip10pt or \]

\[\lim_{n\mapsto +\infty } s(j10_{n}wi)=
\frac{(1+\tau_{1}(a^{*}))^{\frac{1}{\gamma}}-(1+\tau_{2}(a^{*}))^{\frac{1}{\gamma}}}
{(1+\tau_{1}(a^{*}))^{\frac{1}{\gamma}}-1}. \]

\vskip5pt
Because $I_{j10_{n}wi}$ is in one of the middle intervals, $I_{01}$
and $I_{11}$. from the proof of Lemma
A3, the error of $s(j10_{n}i)$ to $s_{0}(b^{*})$ can be estimated
by $|I_{j10_{n}wi}|^{\alpha }$, that is, there is a positive
constant $C_{2}$ such that
\[ |s_{0}(b^{*})
-s(j10_{n}wi)|\leq C_{2} |I_{j10_{n}wi}|^{\alpha }\]
for any $b^{*}=(\cdots j10_{n}wi.)$ in ${\cal C}^{*}$.
Now we can get that the limit of $s_{0}(b^{*})$ as $b^{*}\neq a^{*}$
tends to $a^{*}$ exists and
\[ \lim_{b^{*}\neq a^{*}, b^{*}\mapsto a^{*}}
s_{0}(b^{*})=
\frac{(1+\tau_{2}(a^{*}))^{\frac{1}{\gamma}}-1}
{(1+\tau_{1}(a^{*}))^{\frac{1}{\gamma}}-1} \hskip10pt or \]

\[ \lim_{b^{*}\neq a^{*}, b^{*}\mapsto a^{*}}
s_{0}(b^{*})=
\frac{(1+\tau_{1}(a^{*}))^{\frac{1}{\gamma}}-(1+\tau_{2}(a^{*}))^{\frac{1}{\gamma}}}
{(1+\tau_{1}(a^{*}))^{\frac{1}{\gamma}}-1}.\]

\vskip5pt
Because $s_{0}(a^{*})
=\tau_{2}(a^{*})/\tau_{1}(a^{*})$, the limit of $s_{0}(b^{*})$ as
$b^{*}\neq a^{*}$ tends to $a^{*}$
does not equal $s_{0}(a^{*})$. In other words, $s_{0}$ has jump
discontinuity at $a^{*}$. \Q

\vskip5pt
\nd {\sc Lemma A5.}
{\em The scaling function $s_{0}$ is continuous at
all points in ${\cal B}$.}

\vskip5pt
\nd {\sl Proof.} Suppose $a^{*}$ is in ${\cal B}$.  Let
$b^{*}$ be any point in ${\cal C}^{*}$ with the same first
$(n+1)^{th}$ coordinates $w_{n}i$ as that of $a^{*}$.
If the $n^{th}$ coordinate of $a^{*}$ is one, then $I_{w_{n}}$ is in
one of the middle intervals $I_{01}$ and $I_{11}$.
The errors from $s_{0}(a^{*})$ and
$s_{0}(b^{*})$ to
$s(w_{n}i)$ can be estimated by
$|I_{w_{n}}|^{\alpha }$,  that is,  there is a positive constant
$C$ such that
$|s_{0}(a^{*}) -s(w_{n}i)|$ and
$|s_{0}(b^{*}) -s(w_{n}i)|$ are less than $ C|I_{0,
w_{n}}|^{\alpha }$.
Then $|s_{0}(a^{*})-s_{0}(b^{*})|\leq 2C|I_{w_{n}}|^{\alpha }$.
Because there are infinite many ones in the coordinates of $a^{*}$, the
limit of $\{ s_{0}(b^{*})\} $ as $b^{*}$ tends to $a^{*}$ exists and
$\lim_{b^{*}\mapsto
a^{*}}s_{0}(b^{*}) =s_{0}(a^{*})$. In other words, $s_{0}$ is continuous
at $a^{*}$. \Q

\vskip5pt
Suppose $\tilde{f}_{0}=h_{\gamma}\circ f_{0}\circ h_{\gamma }^{-1}$ is
again the representation of $f_{0}$ under the singular metric
associated to $f_{0}$.

\vskip5pt
\nd {\sc Lemma A6.} {\em There is a H\"older continuous scaling
function $\tilde{s}_{0}$ of $\tilde{f}_{0}$ and
the restriction of $s_{0}$ to
${\cal B}$ equals the restriction of $\tilde{s}_{0}$
to
${\cal B}$. In particular, the restriction of $s_{0}$ to ${\cal B}$ is
H\"older continuous on ${\cal B}$.}

\vskip5pt
\nd {\sl Proof.}
By using similar
arguments
to the proof of Lemma 6,  we can show that
there is a H\"older continuous scaling
function $\tilde{s}_{0}$ of $\tilde{f}_{0}$.

\vskip5pt
The restriction of $h_{\gamma}$ to the union
of the middle intervals $I_{10}$
and $I_{11}$ is a $C^{1}$ embedding.  For any $a^{*}$ in ${\cal B}$, if
the $n^{th}$ coordinate of $a^{*}$ is one, then $I_{
w_{n}}$ is in one of the middle intervals.
We use
$|h_{\gamma }(I_{w_{n}i})|/|h_{\gamma}(I_{w_{n}})| $ and
$|I_{w_{n}i}|/|I_{w_{n}}|$ to approach
$\tilde{s}_{0}(a^{*})$ and $s_{0}(a^{*})$,
respectively.
Now we can show that $\tilde{s}_{0}(a^{*}) =s_{0}(a^{*})$ because there
are infinite many ones in the coordinates of $a^{*}$. \Q

\vskip5pt
Lemma 6 and Lemma A1 to A6 give the proof of
Theorem A.

\vskip10pt
\nd {\bf \S 4.5 Scaling functions of mappings on ${\cal BH}$}

\vskip5pt
More generally,
we
have the following theorem.

\vskip5pt
\nd {\sc Theorem C.} {\em
Suppose $f$ is on ${\cal BH}$ and $\tilde{f}$ is
the representation of $f$ under the singular metric associated to $f$.
There
exist the scaling function $s_{\tilde{f}}$ of
$\tilde{f}$ and the scaling function $s_{f}$ of $f$ and these
scaling functions satisfy that

\vskip5pt
\nd (a) $s_{\tilde{f}}$ is H\"older continuous on ${\cal C}^{*}$,

\vskip5pt
\nd (b) $s_{f}$ has jump discontinuities at
all points in ${\cal A}$ and
$s_{f}$ is continuous at
all points in ${\cal B}$,

\vskip5pt
\nd (c) the restriction of $s_{f}$ to ${\cal B}$ equals the
restriction of $s_{\tilde{f}}$ to ${\cal B}$}.

\vskip5pt
\nd {\sl Proof.} The proof is the same as the proofs of Lemma A3 to
Lemma A6. \Q

\vskip5pt
Suppose ${\cal S}= \{ s_{f} | f\in {\cal BH}\}$.
We can use $s_{f}$ to
compute the eigenvalues of $f$ at all periodic points (see
[J3]) and the power law at the critical point, that is,
\[ \gamma
=\frac{\log s_{f}((0_{\infty }.))}{\log (\lim_{b^{*}\in {\cal B},
b^{*}\mapsto
0_{\infty }} s_{f}(b^{*}))}.\]
The absolute value of the asymmetry of $f$
at the critical point is $|sv_{f}|=\lim_{n\mapsto +\infty
}|I_{010_{n}}|/|I_{110_{n}}|$.

\vskip5pt
An example of a scaling function in ${\cal S}$ is the following
proposition.

\vskip5pt
\nd {\sc Proposition 2.}
{\em  Let $q: x\mapsto 1-2x^{2}$.  Then
$s_{q}(a^{*}) =1/2$ for $a^{*}\in {\cal B}$ and
$s_{q}(a^{*})\neq 1/2$ for $a^{*} \in {\cal A}$.}

\vskip5pt
\nd {\sl Proof.} Recall that $\tilde{q}(y)=1-2|y|$.  The proof of
this proposition just follows the proof of Lemma A3 and Lemma A6. \Q

\vskip10pt
\nd {\bf \S 4.6 The Hausdorff dimension of the maximal invariant set of
$q_{\ep}$.}

\vskip5pt
From Proposition 2, we can observe more on the Hausdorff
dimension of the maximal invariant set of $q_{\ep }(x)=
1+\ep-(2+\ep)x^{2}$.
Suppose $\Lambda_{\ep}$ is the maximal invariant set of $q_{\ep}$ and
$HD(\ep)$ is the Hausdorff dimension of
$\Lambda_{\ep}$.

\vskip5pt
\noindent {\sc Proposition 3.}
{\em
There
is a \poc\ $C$ which does not dependent on $\ep $ such that
\[1- C^{-1}\sqrt{\varepsilon } \leq HD(\ep) \leq
1- C\sqrt{\varepsilon }\]
for all $0\leq \ep \leq 1$.}

\vskip5pt
\nd {\sl Proof.}
Suppose $dy =
dx/\sqrt {(1+\varepsilon )^{2} -x^{2}}$ is the metric associated to
$q_{\ep}$ on
$[-1,1]$,  $y=h_{2,\ep}(x) $ is the corresponding change of
coordinate and $\tilde{q}_{\ep}=h_{2,\ep}\circ q_{\ep}\circ
h^{-1}_{2, \ep }$. We suppress $\ep$ when there can be no
confusion.

\vskip5pt
Let $\ep$ be in $[0,1]$. We call
$I_{0_{4}}$ and $I_{10_{4}}$ the end intervals. Recall that
$I_{0_{4}}$ is the interval in $\eta_{3}$ adjacent to $-1$ and
$I_{10_{4}}$ is the interval in $\eta_{4}$ adjacent to $1$.
We call
the complement of the interiors of the end intervals the middle
interval.

\vskip5pt
Suppose $y$ is in $[-1,1]$ and $x$ is the
preimage of $y$ under $h_{2}$.	The nonlinearity of
$\tilde{q}$ at $y$ is
\[ n(\tilde{q})(y) =
\frac{\ep (1+\ep )}{(2(1+\ep )- (2+\ep)x^{2})\sqrt{(1+\ep )^{2} -x^{2}}}.\]
We can find a positive constant $C_{1}$ which does not depend on
$\ep$ such
that
\[ |n(\tilde{q})(y)| \leq C_{1}\ep\]
for any $y$ in the image of
the middle interval under $h_{2}$.

\vskip5pt
Let $\tilde{\Lambda}$ be the maximal invariant
set of $\tilde{q} $. It is diffeomorphic to $\Lambda$.
The sets
$\tilde{\Lambda}$
and $\Lambda$ have the same Hausdorff dimension. We use $\tilde{I}$
to denote the image of $I$ under $h_{2}$ and computer the
Hausdorff dimension of $\tilde{\Lambda}$.

\vskip5pt
By direct computations,
there is a positive constant $C_{2}$ which does not depend on $\ep$ such
that for any $I_{w}\in \eta_{3}$,
\[ 1-C_{2}^{-1}\ep \leq
|\tilde{I}_{w0}|/|\tilde{I}_{w1}|\leq 1+C_{2}\ep .\]
Suppose $w$ is a finite string of zeroes and ones.
We call a piece string of consecutive zeros
in $w$ a zero element. We call the maximum length of all
the zero elements in
$w$ the zero-length of $w$. If the length of $w$ is greater than $4$ and
$I_{w}$ is in the union of the end intervals,
then the zero-length of $w$ has to be greater then or equal to $4$.
Using this fact,
we can show that for any finite string $w $ of zeros and ones,
if the length of $w$ is greater than $4$ and the zero-length of $w$ is
less than $4$,
the image under the $i^{th}$ iterate of $q$
for any $0< i\leq n-4$
is in the middle interval
for any $0< i\leq n-4$.
For any finite string $w$ of zeroes and ones satisfying that the
length of $w$ is greater then $4$ and
the zero-length of $w$ is less than $3$, the image of $I_{w}$ under
the $(n-4)^{th}$ iterate of $q_{\ep}$ is in $\eta_{3}$ and
$|\tilde{q}_{\ep}^{\circ
(n -4)}(\xi )|/|\tilde{q}_{\ep}^{\circ (n-4)}(\theta )|$ is bounded
above by $\exp(C_{1}\ep )$
for any $\xi$ and $\theta$ in $\tilde{I}_{w}$.
We can find a positive
constant $C_{3}$ which does note depend on $\ep$ such that
\[ 1-C_{3}\ep \leq
\frac{|\tilde{I}_{w0}|}{|\tilde{I}_{w1}|}\leq 1+C_{3}\ep \]
for any finite string $w$ of zeroes and ones satisfying that the
length of $w$ is greater then $4$ and
the zero-length of $w$ is less than $3$.

\vskip5pt
The restriction of the nonlinearity,
$n(h_{2})(x) =x/ ((1+\ep )^{2}-x^{2})$,
to the middle interval is bounded above by a positive constant
$C_{4}$ which does not depend on $\ep$.
For any finite
string $wi$ of zeroes and ones,
we use $\tilde{s}(wi)$ to denote the ratio $|\tilde{I}_{
wi}|/|\tilde{I}_{w}|$.
A direct consequence of Theorem B is
that there is a positive constant $C_{5}$ which does not depend on
$\ep$ such that
$\tilde{s}(w0)
+\tilde{s}
(w1)\geq 1-C_{5}\sqrt{\ep}$ for
all finite string $w$ satisfying that the length of $w$ is
greater then $4$ and the zero-length is less than $3$. Moreover,
there is positive constant $C_{6}$ which does not depend on $\ep$ such
that
$\tilde{s}(wi)\geq (1/2)(1-C_{6}\sqrt{\ep})$
where $i=0$ or $1$.

\vskip5pt
Let $S_{n}=\sum
|\tilde{I}_{w}|^{\delta}$ where sum is over all finite string $w$
satisfying that the zero-length is less than $3$. For
$w=i_{0}i_{1}\cdots i_{n}$, let $w_{k}= i_{0}\cdots i_{k}$ for $4\leq
k\leq n$. We can write $S_{n}$ in
\[ \sum
(\tilde{s}(w_{n})\tilde{s}(
w_{n-1})\cdots \tilde{s}(w_{4}))^{\delta}|\tilde{I}_{
w_{4}}|^{\delta}.\]
Suppose
$C_{7}$ is the minimum length of the intervals in $\eta_{3}$.
Then
$S_{n}$
is greater then
$C_{7}((1/2)(1-C_{6}\sqrt{\ep}))^{\delta(n-4)}
N_{n}$, where $N_{n}$ is the total number of finite strings of zeroes
and ones which satisfies that the zero-length are less than $3$.

\vskip5pt
It is easy to check that $N_{2}=4$
and $N_{n}=2N_{n-1}-1$
for any $n>2$.	We can find a positive constant $C_{8}$
such that $S_{n}$ is greater than $C_{8}((1-C_{6}
\sqrt{\ep})/2)^{n\delta }2^{n}$ for $n\geq 4$.

\vskip5pt
Let $\delta_{0}=\log 2 /(\log 2-\log (1-C_{6}
\sqrt{\ep}))$.	Then $2((1-C_{6}\sqrt{\ep})/2)^{\delta_{0}} =1$ and
$HD(\ep)$ is greater than
or equals to $\delta_{0}$. Here $\delta_{0}$ is bounded below by
$1-sqrt{\ep}$ for a positive constant $C$ which does not
depend on $\ep$.

\vskip5pt
Another side of the inequality in Proposition 3 comes from Corollary 2. \Q

\vskip30pt
\noindent{\Large \bf Acknowledgements}

\vskip5pt
{\em This paper is essentially the second part of my thesis. It owes much 
to the enthusiasm 
and insight of my advisor
Dennis Sullivan. I wish to express my thanks to him, to Wellington de
Melo and Sebastian van Strien and Richard Sacksteder for helpful
conversations
and to Charles Tresser and Frederick Gardiner for reading this paper.
\vskip10pt
I also wish to express my thanks to all members in the dynamical
systems group
at the Graduate Center of the City University of New York.}

\end{document}